\documentclass{article}[12pt]
\usepackage{amsmath,amsfonts,amsthm}
\usepackage{amssymb}

\usepackage{epsfig}

\usepackage{graphics}
\usepackage{graphicx}

\usepackage{latexsym}
\usepackage{graphics}
\usepackage{amsmath}
\usepackage{amsthm}
\usepackage{xspace}
\usepackage{amssymb}
\usepackage{color}
\usepackage{cite}
\usepackage{latexsym}
\usepackage{graphics}
\usepackage{amsmath}
\usepackage{amsthm}
\usepackage{xspace}
\usepackage{amssymb}
\usepackage{epsfig}

\DeclareMathOperator*{\argmax}{arg\,max}
\DeclareMathOperator*{\argmin}{arg\,min}

\newcommand{\beql}[1]{\begin{equation}\label{#1}}
\newcommand{\eeql}{\end{equation}}
\newcommand{\eqn}[1]{(\ref{#1})}

\newcommand{\R}{\mathbb{R}}
\newcommand{\pr}{\mathbb{P}}
\newcommand{\E}{\mathbb{E}}

\newcommand{\ci}{{\cal I}}

\newcommand{\ck}{{\cal K}}
\newcommand{\cx}{{\cal X}}
\newcommand{\cm}{{\cal M}}
\newcommand{\cq}{{\cal Q}}

\newtheorem{thm}{Theorem}
\newtheorem{lem}[thm]{Lemma}

\newtheorem{definition}[thm]{Definition}

\begin{document}

\title{An infinite server system with general packing constraints
%\\(DRAFT)
}

\author
{
Alexander L. Stolyar \\
Bell Labs, Alcatel-Lucent\\
600 Mountain Ave., 2C-322\\
Murray Hill, NJ 07974 \\
\texttt{stolyar@research.bell-labs.com}
}

%\date{June 15, 2009}
\date{\today}

\maketitle

\begin{abstract}

We consider a service system model primarily
motivated by the problem of efficient assignment of virtual machines to physical host machines in 
a network cloud, so that the number of occupied hosts is minimized.

There are multiple input flows of different type customers, with a customer mean service time depending on its type. 
There is infinite number of servers. A server packing {\em configuration} is the vector $k=\{k_i\}$, where $k_i$ is the number 
of type $i$ customers the server "contains". Packing constraints must be observed,
namely there is a fixed finite set of configurations $k$ that are allowed.
Service times of different customers are independent;  
after a service completion, each customer leaves its server and the system.
Each new arriving customer is placed for service immediately; it 
can be placed into a server already serving other customers 
(as long as packing constraints are not violated),
or into an idle server. 

We consider a simple parsimonious real-time algorithm, 
called {\em Greedy}, which attempts to minimize the increment of the objective function 
$\sum_k X_k^{1+\alpha}$, $\alpha>0$, caused by each new assignment; here
$X_k$ is the number of servers in configuration $k$.
(When $\alpha$ is small, $\sum_k X_k^{1+\alpha}$ approximates the total number $\sum_k X_k$ of occupied servers.)
Our main results show that certain versions of 
the Greedy algorithm are {\em asymptotically optimal}, 
in the sense of minimizing $\sum_k X_k^{1+\alpha}$ in stationary regime,
as the input flow rates grow to infinity. 
We also show that
in the special case when the set of allowed configurations  is determined
by {\em vector-packing} constraints, Greedy algorithm can work with
{\em aggregate configurations} as opposed to exact configurations $k$,
thus reducing computational complexity 
while preserving the asymptotic optimality.

\end{abstract}

%\noindent
%{\em Key words and phrases:} Queueing networks, Interacting particle systems,
%Stability, Back-pressure, MaxWeight, Infinite tandem queues, TASEP
%
%\noindent
%{\em Abbreviated Title:} Queue scaling under back-pressure
%
%\noindent
%{\em AMS 2000 Subject Classification:} 
%90B15, 60K25, 60K35, 68M12

%\newpage

\section{Introduction}
\label{sec-intro}

The primary motivation for this work is the following problem arising in cloud computing: 
how to assign various types of virtual machines to physical host machines 
(in a data center) in real time, so that the total number of host machines 
in use is minimized. It is very desirable that an assignment algorithm 
is simple, does need to know the system parameters,
 and makes decisions based on the current system state only.
(An excellent overview of this and other resource allocation issues arising in cloud computing can be found in \cite{Gulati2012}.)

A data center (DC) in the ``cloud'' consists of a number 
%$M$ 
of host machines. 
Assume that all hosts are same: each of them possesses 
the amount $B_n>0$ of resource $n$, where $n\in \{1,2,\ldots,N\}$ is a resource index.
(For example, resource $1$ is CPU, resource $2$ is memory, etc.)
The DC receives requests for virtual machine (VM) placements;  VMs can be of different types $i\in\{1,\ldots,I\}$; a type $i$ VM requires the amounts $b_{i,n}> 0$ of each resource $n$.
Several VMs can share the same host, as long as the host's capacity constraints are not violated; namely, a host can simultaneously contain 
a set of VMs given by a vector $k=(k_1,\ldots,k_I)$, where $k_i$ is the number of type $i$ VMs,
as long as for each resource $n$
\beql{eq-intro1}
\sum_i k_i b_{i,n} \le B_n.
\end{equation}
Thus, VMs can be assigned to hosts already containing other VMs, subject to the above ``packing'' constraints. After a certain random sojourn (service) time each VM vacates its host
(leaves the system), which increases the ``room'' for new arriving VMs to be potentially 
assigned to the host. 
A natural problem is to find a real-time algorithm for assigning VM requests to the hosts, which minimizes (in appropriate sense) the total number of hosts in use. Clearly, such a scheme will maximize the DC capacity; or, if it leaves a large number of hosts unoccupied, those hosts can be (at least temporarily) turned off to save energy. 

More specifically, the model assumptions that we make are as follows:\\
(a) The exact nature of "packing" constraints will not be important -- we just assume that the feasible configuration vectors $k$ (describing feasible sets of VMs that can simultaneously occupy one host) form a finite set $\ck$; and assume monotonicity -- if $k\in\ck$ then so is any $k'\le k$.\\
(b) There is no limit on the number of hosts that can be used and each new VM  is assigned to a host immediately -- so it is an {\em infinite server} model, with no blocking or waiting.\\
(c) Service times of different VMs are independent of each other, even for VMs served simultaneously on the same host.\\
(d)  We further assume in this paper
that the arrival processes of VMs of each type are Poisson
and service time distributions are exponential. 
These assumptions
are not essential and can be much relaxed, as discussed in Section~\ref{sec-more-general}.

The basic problem we address in this paper is:
\beql{eq-problem-basic}
\mbox{minimize} ~~ \sum_k X_k^{1+\alpha},
\end{equation}
where $\alpha> 0$ is a fixed parameter, and 
 $X_k$ is the (random) number of hosts having configuration $k$ in the stationary regime. 
(Clearly, when $\alpha$ is small, $\sum_k X_k^{1+\alpha}$ approximates the total number $\sum_k X_k$ of occupied hosts.)
We consider the {\em Greedy} real-time (on-line) VM assignment algorithm, which,
roughly speaking, tries to minimize the increment of the objective function 
$\sum_k X_k^{1+\alpha}$ caused by each new assignment.
Our main results show that certain versions of 
the Greedy algorithm are {\em asymptotically optimal}, 
as the input flow rates become large or, equivalently, the average number of VMs
in the system becomes large.
%the ratio of the objective under the Greedy algorithm to that under the best possible algorithms, converges to $1$.
We also show (in Section~\ref{sec-aggregation}) that
in the special case when feasible configurations are determined
by constraints \eqn{eq-intro1}, Greedy algorithm can work with
``aggregate configurations'' as opposed to exact configurations $k$,
thus reducing computational complexity 
%the number  of variables that the algorithm needs to maintain can be reduced,
while preserving the asymptotic optimality.

\subsection{Previous work}

Our model is related to the vast literature on the classical {\em stochastic bin packing} problems.
(For a good recent review of one-dimensional bin packing see e.g. \cite{Csirik2006}.)  In particular,
in {\em online} stochastic bin packing problems,
random-size items arrive in the system and need to be placed according to an online algorithm
 into finite size bins; the items never leave or move between bins;
 the typical objective is to minimize the number of occupied bins.
A bin packing problem is {\em multi-dimensional}, when bins and item sizes are vectors; 
the problems with the packing
constraints  \eqn{eq-intro1} are called {\em multi-dimensional
vector packing}  (see e.g. \cite{Bansal2009} for a recent review). 
Bin packing {\em service} systems arise when there is a random in time input flow of random-sized 
items (customers), which need to be served by a bin (server) and leave after a random service time; 
the server can simultaneously process multiple customers as long as 
they can simultaneously fit into it; the customers waiting for service are queued; 
a typical problem is to determine
the maximum throughput under a given (``packing'') algorithm
for assigning  customers for service. 
(See e.g. \cite{Gamarnik2004} for a review of this 
line of work.) Our model is similar to the latter systems, except there are multiple bins (servers), 
in fact -- infinite number in our case.
Models of this type are more recent (see e.g. \cite{Jiang2012,Maguluri2012}).
Paper \cite{Jiang2012} addresses a real-time VM allocation problem, 
which in particular includes packing constraints; 
the approach of \cite{Jiang2012} is close in spirit to Markov Chain algorithms used in 
combinatorial optimization.
%reduces to a Gibbs-sampler-type randomized algorithm. 
Paper \cite{Maguluri2012} is concerned mostly with maximizing throughput of a queueing system 
(where VMs can actually wait for service)
with a finite number of bins.

The asymptotic regime in this paper is such that the input flow rates scale up to infinity.
In this respect, our work is related to the (also vast) literature on queueing systems in the 
{\em many servers} regime. (See e.g. \cite{ST2010_04} for an overview. The name ``many servers''
reflects the fact that the number of servers scales up to infinity as well, 
linearly with the input rates; this condition is irrelevant in our case of infinite number of servers.) 
In particular,
we study {\em fluid limits} of our system, obtained by scaling the  system state down by 
the (large) total number of customers. We note, however, that packing constraints 
are not present in the previous work on the many servers regime, to the best of our knowledge.

\iffalse

The rest of the paper is organized as follows. Section~\ref{sec-model} presents the formal model 
and main result. The ``reduction'' of our problem to the behavior of the infinite system, 
and basic properties of the latter, are given in Section~\ref{sec-basic}. The subcritical and 
supercritical load cases are treated in Sections~\ref{sec-subcritical} and \ref{sec-supercritical},
respectively.  In Section~\ref{sec-gen-input} we remark on more general
input processes. We conclude in Section~\ref{sec-conc}. The appendix contains some of the proofs.

\fi

\section{Model and main results}
\label{sec-model}

We consider a service system with $I$ input customer flows of different types, indexed by $i \in \{1,2,\ldots,I\} \equiv \ci$. 
Each flow $i$ is Poisson with rate $\Lambda_i >0$. Service time of a type $i$ customer is an exponentially distributed random variable with mean $1/\mu_i$.
All input flows and customer service times are mutually independent.
There is an infinite number of servers. Each server can potentially serve more than one customer simultaneously, subject to the following general "packing" constraints. We say that a vector $k = \{k_i, ~i\in \ci\}$ with non-negative integer components
is a server {\em configuration}, if a server can simultaneously serve a combination of different type customers given by vector $k$. The set of all configurations is finite, and is denoted by $\bar\ck$.
We assume that $k\in \bar\ck$ implies that all "smaller" configurations $k'\le k$ belong to $\bar\ck $ as well. Without loss of generality assume that $e_i \in \bar\ck$ for all types $i$, where $e_i$ is the $i$-th coordinate unit vector (otherwise, $i$-customers cannot be served at all).
By convention, the (component-wise)  zero vector $k=0$ belongs to $\bar\ck$ -- this is the configuration of an "empty" server; we denote by $\ck=\bar\ck \setminus \{0\}$ 
the set of configurations {\em not} including the zero configuration. 
\iffalse
COM: CHECK IF NEED THIS
Also without loss of generality, we assume that all configurations 
$k\in\ck$ ``communicate'' with each other, in the sense that there is a path from
any $k$ to any $k'$, along the lattice edges, that stays within $\ck$ (in particular, it does
{\em not} visiting $0$); otherwise, the set customer types $\ci$ is broken down into subsets
such that no two customers from different subsets can be served simultaneously, and so the model decomposes into independent models for each of the subsets separately.
\fi

An important feature of the model is that simultaneous service does {\em not} affect the service rates of individual customers; in other words, the service time of a customer is unaffected by whether or not there are other customers served simultaneously by the same server. Each arriving customer is immediately placed for service in one of the servers; it can be "added" to an empty or non-empty server as long as configuration feasibility constraint is not violated,
i.e. it can be added to any server whose configuration $k\in\bar\ck$ (before the addition) is such that $k+e_i \in \ck$. 
When the service of an $i$-customer by the server in configuration $k$ is completed,
the customer leaves the system and the server's configuration changes to $k-e_i$.
Denote by $X_k$ the number of servers with configuration $k\in \ck$. The system state is then the vector $X = \{X_k, ~k\in \ck\}$. 

A {\em service discipline} ("packing rule") determines where an arriving customer is placed, as a function of the current system state $X$. Under any well-defined service discipline, the system state at time $t$, $X(t)$, is a continuous time, countable Markov chain. It is easily seen to be irreducible and positive recurrent; the positive recurrence follows from the fact
that the total number $Y_i(t)$ of type $i$ customers in the system is the process independent of the service discipline, and its stationary distribution is Poisson with mean $\Lambda_i/\mu_i$. Therefore, the process $X(t), ~t\ge 0$, has a unique stationary distribution.

We are interested in finding a service disciplines minimizing (in a certain sense) 
the total number of non-empty servers in the stationary regime. For example,
an objective can be
\beql{eq-objective-mean}
\min \E \sum_{k\in\ck} X_k^{1+\alpha}(\infty), 
\end{equation}
where $\alpha\ge 0$ is a parameter, and $X(\infty)$ denotes the random system state in stationary regime. Another possible objective is
\beql{eq-objective-prob}
\min \pr \{ \sum_{k\in\ck} X_k^{1+\alpha}(\infty) > C\}, 
\end{equation}
where $C$ is a fixed threshold. In both problems, setting $\alpha=0$ 
obviously corresponds to the exact objective of minimizing the 
number of servers in use; however, if a good discipline for $\alpha>0$ exists,
using such discipline with positive $\alpha$ close to $0$ would (hopefully)
also produce a good solution for the $\alpha=0$ case.

In this paper we consider objectives \eqn{eq-objective-mean} and \eqn{eq-objective-prob}
with $\alpha>0$, and prove the following {\em Greedy} service discipline
(or, rather its slight modification, to be precise) is 
asymptotically optimal as the flows' input rates become large.

\begin{definition}[Greedy discipline]  $ $ 
\begin{enumerate}

\item {\em Integral form (Greedy-I)}. Define $F(x)=\sum_{k\in\ck} (1+\alpha)^{-1} x_k^{1+\alpha}$, $x\in \R_+^{|\ck|}$. A type $i$ customer arriving at time $t$ is 
added to an available configuration $k$ (with either $k=0$ or $X_k(t-)>0$)
such that $k+e_i\in \ck$ and
the increment $F(X(t))-F(X(t-))$ is the smallest.
(Here $X(t-)$ and $X(t)$ are the states just before and just after the addition;
so that $X_{k+e_i}(t)=X_{k+e_i}(t-)+1$ and, unless $k=0$, $X_k(t)=X_k(t-)-1$.)
The ties are broken according to an arbitrary deterministic rule.

\item {\em Differential form (Greedy-D)}. For each $k\in\ck$,
denote
$W_k(x)=(\partial / \partial x_x) F(x)=x_k^{\alpha}, ~x\in\R_+^{|\ck|}$. 
A type $i$ customer arriving at time $t$ is 
added to an available configuration $k$ (with either $k=0$ or $X_k(t-)>0$) 
such that $k+e_i\in \ck$ and the difference
$W_{k+e_i}(X(t-))-I\{k\ne 0\}W_k(X(t-))$ is the smallest.
(Here $I\{\cdot\}$ is the indicator function, equal to $1$ when the condition holds and $0$ otherwise.)
The ties are broken according to an arbitrary deterministic rule.

\end{enumerate}

\end{definition}

In the asymptotic analysis of this paper, the two forms 
of Greedy algorithm are essentially identical, because they induce the same
dynamics of the system in the ``fluid limit''. 
We will analyze the differential form, 
as it is slightly more convenient to work with, 
and is probably more easily implementable in practice; it should be clear that 
all results (along with essentially same proofs) hold for the integral form as well.

{\bf Remark.} All results will hold for a more general objective function
$F(x)=\sum_{k\in\ck} c_k x_k^{1+\alpha}$, with arbitrarily positive weights $c_k$.
The generalization is completely straightforward -- we choose to work with 
$c_k=(1+\alpha)^{-1}$ simply to avoid ``carrying'' factors $c_k(1+\alpha)$, 
which clog notation.

We now define the asymptotic regime. Let $r\to\infty$ be a positive scaling parameter.
(To be specific, assume that $r\ge 1$, and $r$ increases to infinity along a discrete
sequence.)
Input rates scale linearly with $r$; namely,
for each $r$, $\Lambda_i = \lambda_i r$, where $\lambda_i$ are positive parameters.
Let $X^r(\cdot)$ be the process associated with system with parameter $r$, and 
$X^r(\infty)$ be the (random) system state in the stationary regime.
For each $i$ denote by $Y^r_i(t) \equiv \sum_{k\in\ck} k_i X^r_k(t)$ the total number
of customers of type $i$. Since arriving customers are taken for service immediately
and their service times are independent (of the rest of the system), the distribution
of $Y^r_i(\infty)$ is Poisson with mean $r \rho_i$, where $\rho_i\equiv \lambda_i/\mu_i$.
Moreover, $Y^r_i(\infty)$ are independent across $i$.
Since the total number of occupied servers is no greater than the total number 
of 
customers, 
$\sum_k  X_k^r(t) \le Z^r(t)\equiv \sum_i Y^r_i(t)$, we have a simple upper bound 
on the total number of occupied servers in steady state,
$\sum_k  X_k^r(\infty) \le Z^r(\infty)$,
% this means that $\sum_k  X_k^r(\infty)$
%is stochastically upper bounded by $\Pi(r \sum_i \rho_i)$, 
where $Z^r(\infty)$ is a Poisson random variable with mean $r \sum_i \rho_i$.
Without loss of generality, from now on in the paper we assume $\sum_i \rho_i=1$.
(This is equivalent to rechoosing parameter $r$ to be $r \sum_i \rho_i$.)

The {\em fluid scaled} process is $x^r(t)=X^r(t)/r$; 
for any $r$, $x^r(t)$ takes values in (the positive orthant of) Euclidean space $\R^{|\ck|}$,
where $|\ck|$ is the cardinality of $\ck$.
Similarly, $y^r_i(t)=Y^r_i(t)/r$ and  $z^r(t)=Z^r(t)/r$.
Since $\sum_k  x_k^r(\infty) \le z^r(\infty)=Z^r(\infty)/r$,
we see that the random variables
$(\sum_k x_k^r(\infty))^{1+\alpha}$ are uniformly integrable in $r$ 
(for any fixed $\alpha\ge 0$).
This in particular implies that the sequence of distributions of $x^r(\infty)$ is tight,
and therefore there always exists a limit in distribution $x^r(\infty)\implies x(\infty)$,
along a subsequence of $r$. 
(The limit depends on the service discipline, of course.)
The limit (random) vector $x(\infty)$ satisfies the following conservation laws:
\beql{eq-cons-laws}
\sum_{k\in\ck} k_i x_k(\infty) \equiv y_i(\infty) = \rho_i, ~~\forall i,
\end{equation}
implying, in particular,
\beql{eq-cons-laws2}
z_i(\infty)\equiv \sum_i y_i(\infty) \equiv \sum_i \rho_i
%\equiv \sum_i \sum_{k\in\ck} k_i x_k(\infty) 
= 1.
\end{equation}
Therefore, the values of $x(\infty)$ are confined to the convex compact 
$|\ck|-I$-dimensional 
%(with $|\ck|$ being the cardinality of $\ck$) 
polyhedron
$$
\cx \equiv \{x\in \R^{|\ck|} ~|~ x_k\ge 0, ~\forall k\in\ck; 
~ \sum_k k_i x_k = \rho_i, ~\forall i\in\ci \}.
%\cx \equiv \R_+^{|\ck|} \cap \left[ \cap_i \{x~|~\sum_k k_i x_k = \rho_i\}\right].
$$
(We will slightly abuse notation by using symbol $x$ for a generic element of $\cx$;
while $x(\infty)$, and later $x(t)$, refer to random variables taking values in $\cx$.)

The asymptotic regime and the associated basic properties \eqn{eq-cons-laws}
and \eqn{eq-cons-laws2} hold {\em for any service discipline}, not necessarily
Greedy-D.

%Let us consider $F(x)=\sum_k (1+\alpha)^{-1} x_k^{1+\alpha}$, where $\alpha>0$ is fixed. 
Note that function
$F(x)$ with $\alpha>0$ is strictly
convex on $\cx$, and therefore there is the unique optimal point $x^*$ minimizing $F$:
\beql{eq-xstar}
x^* = \argmin_{u\in\cx} F(u).
\end{equation}
(Note that, if $\alpha=0$, then $x^*$ is an optimal solution of a linear program and therefore might not be unique.
We do not consider the $\alpha=0$ case in this paper.)

Our {\bf main results} are as  follows.\\
 {\em Let $\alpha>0$. We prove that, as $r\to\infty$, the convergence 
is distribution $x^r(\infty) \implies x^*$ holds in two cases: \\
(a) For the closed system, with the constant number $Y^r_i=\rho_i r$ of customers of each type, operating under the Greedy-D discipline. 
(The exact result is Theorem~\ref{th-fluid-stat-gt1-closed}.)\\
(b) For the original system (as defined above), operating under a slightly modified
Greedy-D discipline, called Greedy-DM. 
(The exact result is Theorem~\ref{th-fluid-stat-gt1-Greedy-DM}.)\\
In addition, in the special case when feasible configurations are determined
by vector-packing constraints \eqn{eq-intro1}, 
we prove that essentially same results hold if Greedy algorithm
uses quantities $X^r_k$ aggregated over classes of equivalent configurations,
thus reducing the total number of system variables the algorithm needs to maintain.
(The exact results are Theorems~\ref{th-fluid-stat-gt1-closed-aaa}
and \ref{th-fluid-stat-gt1-Greedy-DM-aaa}.)
}

\iffalse
\begin{thm}
\label{th-fluid-stat-gt1}
Suppose $\alpha>0$ is fixed, and
for each $r$ Greedy-D with 
$W_k(X)=X_k^\alpha$, is applied. Then, as $r\to\infty$,
$$
x^r(\infty) \implies x^*.
$$
\end{thm}
\fi

\subsection{Basic notation and conventions}
\label{subsec-notation}

Standard Euclidean norm of a vector $x\in \R^n$ is denoted $\|x\|$; 
the distance from vector $x$ to a set $U$ in a Euclidean space is denoted
$d(x,U)=\inf_{u\in U} \|x-u\|$; $\R_+$ is the set of real non-negative numbers;
the cardinality of a finite set $\mathcal{N}$ is 
$|\mathcal{N}|$.
Symbol $\to$ means ordinary convergence in $\R^n$; $\implies$
denotes convergence in distribution of random variables taking values in $\R^n$,
equipped with the Borel $\sigma$-algebra;
abbreviation {\em w.p.1} means convergence {\em with probability 1}.
We often write $x(\cdot)$ to mean the function (or random process) $x(t),~t\ge 0$.
Abbreviation {\em u.o.c.} means 
{\em uniform on compact sets} convergence of functions, with the argument (usually in 
$[0,\infty)$) determined by the context.
We often  write $\{x_k\}$ to mean the vector $\{x_k, ~k\in\ck\}$, with the set of indices
$\ck$ determined by the context.

For a finite set of scalar functions $f_n(t), ~t\ge 0$, $n\in\mathcal{N}$, a point $t$ is called
{\em regular} if for any subset $\mathcal{N}' \subseteq \mathcal{N}$ the proper
derivatives
$$
\frac{d}{dt} \max_{n\in\mathcal{N}'} f_n(t) ~~\mbox{and}~~ 
\frac{d}{dt} \min_{n\in\mathcal{N}'} f_n(t)
$$
exist.

\section{Closed system. Greedy-D optimality
%Proof of Theorem~\ref{th-fluid-stat-gt1}
}
\label{sec-proof-closed-gt1}

In this section we consider a ``closed'' version of our system.
Namely, assume that there is a fixed number $\rho_i r$ customers of type $i$
(in a system with parameter $r$); there are no exogenous arrivals into the system --
when a service of a type $i$ customer is completed, the customer
 immediately has to be placed 
into a server for a new service. Service discipline determines where the customer 
is placed, based on the current system state.

It is easy to see that, for any $r$, a stationary distribution of the process exists
under any given service discipline,
because the process in this case is a finite-state continuous time Markov chain.

The main result of this section is the following 
%theorem, showing that 
%when $r$ is large, $x^r(t)$ is steady-state is close to $x^*$.

\begin{thm}
\label{th-fluid-stat-gt1-closed}
Consider a sequence of closed systems, indexed by $r$,
and let $x^r(\infty)$ denote the random state of the (fluid-scaled) process
is a stationary regime, under the Greedy-D discipline with $\alpha>0$.
%with  $W_k(X)=X_k^\alpha$, $\alpha>0$.
Then, as $r\to\infty$,
$$
x^r(\infty) \implies x^*,
$$
where $x^*$ is defined in \eqn{eq-xstar}.
\end{thm}

To prove the theorem we will need to study the transient behavior 
of the fluid-scaled process and its limits.

Let $\cm$ denote the set of pairs $(k,i)$ such that $k\in\ck$ and $k-e_i\in\bar\ck$.
Each pair $(k,i)$ is associated with the ``edge'' $(k-e_i,k)$ connecting configurations
$k-e_i$ and $k$; often we refer to this edge as $(k,i)$.
By ``arrival along the edge $(k,i)$'' we will mean placement of a type $i$ customer
into a server configuration $k-e_i$ to form configuration $k$;
similarly, ``departure along the edge $(k,i)$'' is a departure of a type $i$ customer
from a server in configuration $k$, which changes its configuration to $k-e_i$.

For each $(k,i)\in \cm$, consider an independent
unit-rate Poisson process $\Pi_{ki}(t), ~t\ge 0$.
We have the functional strong law of large numbers:
\beql{eq-flln-poisson}
\frac{1}{r}\Pi_{ki}(rt) \to t, ~~~u.o.c., ~~w.p.1.
\end{equation}

Without loss of generality, assume that
the Markov process $X^r(\cdot)$ for each $r$ is driven by the common
set of Poisson processes $\Pi_{ki}(\cdot)$, as follows.
For each $(k,i)\in \cm$, let us denote by $D^r_{ki}(t)$ the total 
number of  departures along the edge $(k,i)$ in $[0,t]$; then we
can assume that
\beql{eq-driving}
D^r_{ki}(t) = \Pi_{ki} (\int_0^t X_k^r(\xi) k_i \mu_i d\xi).
\end{equation}
Each type-i departure in the closed system is simultaneously a type-i ``arrival'',
which is allocated according to Greedy-D. Thus, the realization of the process
is (w.p.1) uniquely determined by the initial state $X^r(0)$ and the realizations 
of $\Pi_{ki}(\cdot)$. Denote by $A^r_{ki}(t)$ the total number of arrivals 
allocated along edge $(k,i)$. Obviously, we have the conservation law for each type $i$:
$\sum_k A^r_{ki}(t) \equiv \sum_k D^r_{ki}(t), ~\forall t\ge 0$.
In addition to
$$
x^r_{k}(t) = \frac{1}{r} X^r_{k}(t),
$$
we introduce other  fluid-scaled 
quantities:
$$
d^r_{ki}(t) = \frac{1}{r} D^r_{ki}(t),~~~
a^r_{ki}(t) = \frac{1}{r} A^r_{ki}(t).
$$

A set of Lipschitz continuous functions
$[\{x_k(\cdot),~k\in \ck\}, \{d_{ki}(\cdot),~(k,i)\in \cm\},\{a_{ki}(\cdot),~(k,i)\in \cm\}]$
on the time interval $[0,\infty)$ we call a {\em fluid sample path} (FSP), if there exist
realizations of $\Pi_{ki}(\cdot)$ satisfying \eqn{eq-flln-poisson}
and a fixed subsequence of $r$, along which
%\beql{eq-fsp-def-closed}
\begin{eqnarray}
& [\{x_k^r(\cdot),~k\in \ck\}, \{d_{ki}^r(\cdot),~(k,i)\in \cm\},\{a_{ki}^r(\cdot),~(k,i)\in \cm\}]
\to \nonumber \\
& [\{x_k(\cdot),~k\in \ck\}, \{d_{ki}(\cdot),~(k,i)\in \cm\},\{a_{ki}(\cdot),~(k,i)\in \cm\}],
~~u.o.c. \label{eq-fsp-def-closed}
\end{eqnarray}
%\end{equation}

It is easy to see that the family of all FSPs is {\em uniformly} Lipschitz.

\begin{lem}
\label{lem-conv-to-fsp-closed}
Suppose the initial states $x^r(0)$ are fixed and are such that $x^r(0)\to x(0)$.
Then, w.p.1 for any subsequence of $r$ there exists a further
subsequence of $r$, along which the convergence \eqn{eq-fsp-def-closed} holds,
where the limit is an FSP.
\end{lem}

{\em Proof} is very standard.
Essentially, it suffices to observe that, with probability 1, each sequence (in $r$)
of functions $d^r_{ki}(\cdot)$ is asymptotically Lipschitz, namely for some $C>0$, 
and all $0\le t_1 \le t_2 <\infty$, 
$$
\limsup_r d^r_{ki}(t_2) - d^r_{ki}(t_1) \le C(t_2-t_1),
$$
which in turn follows from \eqn{eq-flln-poisson}. 
And similarly for functions $a^r_{ki}(\cdot)$, because $a^r_{ki}(t) \le \sum_k d^r_{ki}(t)$.
Using this and \eqn{eq-fsp-def-closed}, we
easily verify the u.o.c. convergence of all $x^r_k(\cdot),d^r_{ki}(\cdot),a^r_{ki}(\cdot)$,
along possibly a further subsequence,
and the fact that the limits are Lipschitz.
We omit details.
$\Box$

\iffalse
COM: Essentially, it suffices to observe that, since each $x^r_k(t)\le 1$,
all functions $\int_0^t x^r_k(\xi)d\xi$, for all $r$, are uniformly Lipschitz,
and therefore there always exits a subsequence of $r$ along which they converge u.o.c.
to Lipschitz functions. Using this, along with \eqn{eq-fsp-def-closed}, we
easily verify the convergence of all $x^r_k(\cdot),d^r_{ki}(\cdot),a^r_{ki}(\cdot)$,
along possibly a further subsequence,
and the fact that the limits are Lipschitz.
\fi

For an FSP, at a regular time point $t$, we denote
$v_{ki}(t)=(d/dt)a_{ki}(t)$ and $w_{ki}(t)=(d/dt)d_{ki}(t)$.
In other words, $v_{ki}(t)$ and $w_{ki}(t)$ are the rates of type $i$ ``fluid''
arrival and departure along edge $(k,i)$, respectively.

\begin{lem}
\label{lem-fsp-properties-closed-basic}
An FSP satisfies the following properties:
$$
y_i(t)=\sum_k k_i x_k(t) \equiv \rho_i;
$$
at any regular point $t$,
$$
w_{ki}(t)=k_i \mu_i x_k, ~~\forall (k,i)\in \cm,
$$
$$
\sum_{k:(k,i)\in \cm} w_{ki}(t) = \sum_{k:(k,i)\in \cm} v_{ki}(t) = \lambda_i, ~~\forall i\in \ci,
$$
$$
(d/dt) x_k(t) = [\sum_{i:k-e_i\in\bar\ck} v_{ki} - \sum_{i:k+e_i\in\ck} v_{k+e_i,i} ]
                     - [\sum_{i:k-e_i\in\bar\ck} w_{ki} - \sum_{i:k+e_i\in\ck} w_{k+e_i,i} ], ~~\forall k\in\ck.
$$
\end{lem}

{\em Proof} is both standard and obvious.
$\Box$

\section{Characterization of the optimal point $x^*$ 
%in \eqn{eq-xstar}, 
and related properties}
\label{sec-opt-greater1}

This section describes properties of the optimal point $x^*$, 
and related general properties 
of the ``allocation'' vectors,
which, roughly speaking, have the meaning of the vector $v(t)= \{v_{ki}(t),~(k,i)\in\cm\}$
of arrival rates. 
The results of this section are {\em not} limited the closed system or 
the Greedy-D algorithm.

Recall that $x^*$ is defined as the unique optimal solution of the convex optimization
problem 
\beql{eq-opt-greater1}
\min_{x\in\R_+^{|\ck|}} F(x)
\end{equation}
subject to
\beql{eq-cons-laws222}
\sum_{k\in\ck} k_i x_k = \rho_i, ~~\forall i,
\end{equation}
%\beql{eq-nonneg222}
%x_k \ge 0, ~~\forall k\in \ck,
%\end{equation}
where $F(x)=\sum_k (1+\alpha)^{-1} x_k^{1+\alpha}$ with $\alpha>0$.

Using Lagrange multipliers $\eta_i$ for the constraints \eqn{eq-cons-laws222},
the Lagrangian is
$$
\sum_k \frac{1}{1+\alpha} x_k^{1+\alpha} + \sum_i \eta_i [\rho_i - \sum_k k_i x_k]
= \sum_k [\frac{1}{1+\alpha} x_k^{1+\alpha} - x_k\sum_i \eta_i k_i] + \sum_i \eta_i \rho_i.
$$
Therefore, we have the following characterization: vector $x=x^*$ 
%is the optimal solution of
%\eqn{eq-opt-greater1}-\eqn{eq-cons-laws222} 
if and only if there exist constants
$\eta_i$ such that
\beql{x-star-characterize}
x_k^{\alpha} = \max\{\sum_i k_i \eta_i,0\}.
\end{equation}
From here we also observe that at least one of the Lagrange multipliers is strictly positive,
$\eta_i>0$, and therefore, for such $i$, $x^*_{e_i}>0$.

For $x \in \R_+^{|\ck|}$ and $(k,i)\in \cm$ denote
$$
\Delta_{ki} = \Delta_{ki}(x) = x_k^{\alpha} - x_{k-e_i}^{\alpha},
$$
where, by convention, $x_k=0$ when $k=0$. 

\begin{definition}[Allocations]
A vector
% (``exogenous input rates'') 
$\gamma=\{\gamma_{ki}, ~ (k,i)\in\cm\}$, such that all $\gamma_{ki}\ge 0$,
$$
\sum_{k ~:~ (k,i)\in \cm} \gamma_{ki} = \lambda_i, ~~\forall i,
$$
we will call an {\em allocation} (of the arrival rates).
%(Component $\gamma_{ki}$ has the meaning of exogenous arrival rate of type $i$ fluid, 
%allocated
%to edge $(k,i)$, that is the rate at which new $i$-arrivals are added to configurations
%$k-e_i$ to form configurations $k$.)
For a given $x\in \cx$,
the allocation $\gamma'=\gamma'(x)$, with components $\gamma'_{ki} = k_i \mu_i x_k$,
is called {\em neutral}. 
An allocation $\gamma$ is called a {\em simple improving allocation},
or {\em SI-allocation},
for a given $x\in \cx$, if 
there exist edges $(k,i)$ and $(k',i)$, with the same $i$ but $k'\ne k$, 
%and number $\delta> 0$
such that the following conditions hold:
$$
k_i>0, ~ x_k>0, %~\delta \le k_i \mu_i x_k,
$$
\beql{eq-key-for-realoc}
\mbox{either} ~~ k'=e_i ~~\mbox{or} ~~ [k'_i>0 ~\mbox{and}~x_{k'-e_i}>0],
\end{equation}
\beql{eq-Delta-ineq}
\Delta_{k'i} < \Delta_{ki},
\end{equation}
$$
\gamma_{ki}= 0, ~~ \gamma_{k'i}= k'_i \mu_i x_{k'} + k_i \mu_i x_k,
$$
$$
\gamma_{\ell j}= \gamma'_{\ell j}= \ell_j \mu_j x_{\ell}, ~~\mbox{for}~(\ell,j)\ne (k,i),(k',i).
$$
We denote by $\Gamma(x)$
the set consisting of the neutral and all SI-allocations for $x$.
\end{definition}

The meaning of this definition is simple and becomes clear with the use of
following notation.
For $x \in \cx$,
denote
$$
D(\gamma)=D(\gamma,x)=\sum_{(k,i)\in\cm} \Delta_{ki} [\gamma_{ki} - k_i \mu_i x_k].
$$
(The meaning of this is: $D(v(t),x(t))=(d/dt)F(x(t))$
for an FSP at a regular point $t$.)
Clearly, $D(\gamma')=0$ for the neutral allocation, 
and $D(\gamma)<0$ for any SI-allocation.
This is because any SI-allocation $\gamma$, associated with edges $(k,i)$ and $(k',i)$,
is obtained from the neutral by ``reallocating'' the positive
amount $k_i \mu_i x_k$ of type-$i$ input
rate from edge $(k,i)$ to edge $(k',i)$ with strictly smaller $\Delta_{k'i}$;
condition \eqn{eq-key-for-realoc} guarantees that there are  
servers in the configuration $k'-e_i$, to which these reallocated type-$i$ arrivals
can be added.

\begin{lem}
\label{lem-F-decrease-key}
If $x\in \cx$ and $x\ne x^*$, then there exists  
at least one SI-allocation $\gamma\in\Gamma(x)$.
\end{lem}

%COM: The proof does use that fact that $\partial F /\partial x_k =0$ when $x_k=0$.

{\em Proof} is by contradiction. Suppose SI-allocations do not exist.
Then there must exist at least one $i\in\ci$ such that $x_{e_i}>0$. If not,
we pick a minimal $k$ with $x_k>0$, for which necessarily $\sum_i k_i\ge 2$.
Then we could pick $i$ such that $k_i\ge 1$ and
construct the SI-allocation $\gamma$ 
associated with edges $(k,i)$ and $(e_i,i)$.
(I.e., $\gamma$ is obtained from the neutral allocation by ``reallocating''
amount $k_i \mu_i x_k$ of type $i$ input from edge $(k,i)$ to edge $(e_i,i)$.)
Therefore, indeed, $x_{e_i}>0$ for at least one $i$.

Denote by $\ci^+$ the set of those $i$ with $x_{e_i}>0$, 
and by $\ci^- = \ci \setminus \ci^+$ the set of those $i$ with $x_{e_i}=0$.
Set $\ci^+$ is non-empty, while $\ci^-$ may or may not be empty.
Let us fix the following values of $\eta_i$:
$$
\eta_i= \left\{ \begin{array}{ll}
               x_{e_i}^{\alpha}, & \mbox{if}~i\in\ci^+\\
               \min_{k~:~(k,i)\in\cm} \Delta_{ki}, & \mbox{if}~i\in\ci^-
                \end{array}
       \right.
$$
It follows from this definition that $\eta_i>0$ if and only if $i\in\ci^+$.

Let us show that \eqn{x-star-characterize} must hold for these $\eta_i$.
Denote $u_k=\max\{\sum_i k_i \eta_i,0\}$; so, we will show that $x_k^{\alpha}=u_k$ for all $k$.
Suppose not. Let us choose a minimal $k$ for which this relation fails.
Note that necessarily $\sum_i k_i\ge 2$.

Consider first the case when $k_i=0$ for all $i\in\ci^-$. 
Pick any $i\in\ci^+$ for which $k_i\ge 1$.
If $x_k^{\alpha}>u_k$, then $\Delta_{ki}>\eta_i$; we can then construct the SI-allocation
associated with the reallocation from $(k,i)$ to $(e_i,i)$;
this is a contradiction. If $x_k^{\alpha}<u_k$, then $\Delta_{ki}<\eta_i$,
and we can construct the SI-allocation, which does the ``opposite''
reallocation from $(e_i,i)$ to $(k,i)$; again, a contradiction.

The second case is when $k_i\ge 1$ for at least one $i\in\ci^-$. We pick such an $i$.
Condition $x_k^{\alpha}<u_k$ would imply
$\Delta_{ki}<\eta_i$ -- a contradiction with the definition of $\eta_i$.
Therefore, $x_k^{\alpha}>u_k$, and then
 $x_k>0$ and $\Delta_{ki}>\eta_i$. If $\eta_i=0$,
then we can construct the SI-allocation 
associated with the reallocation from $(k,i)$ to $(e_i,i)$.
Otherwise, if $\eta_i<0$, there exists an edge $(k',i)$ such that 
$\Delta_{k'i}=\eta_i<0$ and, consequently, $x_{k'-e_i}>0$;
we can construct the SI-allocation 
associated with the reallocation from $(k,i)$ to $(k',i)$.

We have proved that, indeed, \eqn{x-star-characterize} holds for the
chosen values of $\eta_i$. But this means that $x=x^*$ 
-- this contradiction completes the proof. 
$\Box$

For any $x$, denote
$$
D_{min}(x) = \min_{\gamma\in\Gamma(x)} D(\gamma,x).
$$
The minimum is attained and, obviously, $D_{min}(x)\le 0$ for any $x$.

\begin{lem}
\label{lem-F-decrease-key2}
For any $\epsilon>0$ there exists $\epsilon_1>0$ such that
\beql{eq-decay-rate-upper-bound}
\|x-x^*\| \ge \epsilon ~~\mbox{implies}~~ D_{min}(x) < -\epsilon_1.
\end{equation}
\end{lem}

{\em Proof.} We use compactness of $\cx$. If the lemma statement does not hold, 
then, for some fixed $\epsilon>0$ we can find a converging sequence $x^{(n)}\to x'\in\cx$,
$n\to\infty$,
such that $D_{min}(x^{(n)})\to 0$ and $\|x-x'\| \ge \epsilon$. 
%Denote $\ck^+ = \{k ~|~ x'_k>0\}$. 
%Clearly, for all sufficiently large $n$, $x_k^{(n)}>0$ when $k\in\ck^+$.
By Lemma~\ref{lem-F-decrease-key}, there exists an SI-allocation $\gamma\in\Gamma(x')$.
Since $x'_k>0$ implies $x_k^{(n)}>0$ for all large $n$,
we can construct (for all large $n$) an SI-allocation $\gamma^{(n)}\in\Gamma(x^{(n)})$ 
associated with the same reallocation as the one producing $\gamma$.
Then, we have $\gamma^{(n)}\to \gamma$ and $D(\gamma^{(n)},x^{(n)})\to D(\gamma, x')<0$.
This contradicts the assumption $D_{min}(x^{(n)})\to 0$.
$\Box$

\section{Proof of Theorem~\ref{th-fluid-stat-gt1-closed}}
\label{sec-proof-closed}

\begin{lem}
\label{lem-F-decrease}
Any FSP (for the closed system under Greedy-D) is such that
at any regular point $t$, 
%the actual allocation $v=v(t)$ is such that
\beql{eq-F-decrease-closed}
\frac{d}{dt}F(x(t)) = D(v(t),x(t)) \le D_{min}(x(t)).
\end{equation}
\end{lem}

{\em Proof.} Within this proof, $x=x(t)$ and $v=v(t)$.
Consider a specific allocation $\gamma\in\Gamma(x)$, 
for which the minimum in the definition
of $D_{min}(x)$ is attained; unless $x=x^*$, $\gamma$ is an SI-allocation,
associated with some fixed reallocation from $(k,i)$ to  $(k',i)$.
Define the following ``distribution function'':
$$
H(\xi;v) = \sum_{(\ell,j)\in \cm ~:~ \Delta_{\ell j} \le \xi} v_{\ell j}, ~~~\xi\in \R.
$$
Function $H(\xi;\gamma)$ is defined the same way, by replacing $v$ with $\gamma$.
Then, we must have
\beql{eq-H-order}
H(\xi;v) \ge H(\xi;\gamma), ~~\forall \xi.
\end{equation}
Indeed, consider any two edges $(\ell,j)\ne(\ell',j)$ with a common $j$,
and suppose $\Delta_{\ell' j} < \Delta_{\ell j}$. Then, in a fixed sufficiently small
time interval $[t,t+\epsilon]$,
for all pre-limit trajectories with sufficiently large parameter $r$,
any $j$-customer whose service completes at a server with configuration $\ell'$
cannot possibly be placed along the edge $(\ell,j)$.
%i.e. added into configuration $\ell'-e_j$ thus forming $\ell'$.
(Here we use the fact that the system is closed: when a customer service is completed,
there is always an option to place it back into the same server for the new service.)
Furthermore, if either $\ell'=e_j$ or $x_{\ell'-e_j}>0$, any $j$-customer completing service
in configuration $\ell$, will be placed for new service either along the
edge $(\ell',j)$ or possibly along another edge $(\ell'',j)$ with 
$\Delta_{\ell'' j} \le \Delta_{\ell' j}$.
\iffalse
We see that, definitely, $H(\xi;v) \ge H(\xi;\gamma')$.
% i.e. allocation $v$ is at least as ``good'' as neutral. 
(This covers the case $x=x^*$.) Furthermore, when $x\ne x^*$,
$v$ is such that, compared to $\gamma'$,
 at least the ``amount'' $k_i \mu_i x_k$ of type $i$ input rate 
will be reallocated
from edge $(k,i)$  to the edge $(k',i)$ with smaller difference
$\Delta_{k'i}$, or perhaps to other edges $(k'',i)$ with differences 
$\Delta_{k'' i}\le \Delta_{ki}$. 
\fi
This proves \eqn{eq-H-order}, from which \eqn{eq-F-decrease-closed}
easily follows.
$\Box$

As a corollary of Lemma~\ref{lem-F-decrease} we obtain

\begin{lem}
\label{lem-fluid-convergence-closed} 
Any FSP (for the closed system under Greedy-D) is such that
\beql{eq-fluid-convergence-closed}
x(t)\to x^*.
\end{equation}
The convergence is uniform across all initial states $x(0)\in\cx$. 
\end{lem}

{\em Remainder of the proof of Theorem~\ref{th-fluid-stat-gt1-closed}.} 
We fix $\epsilon>0$ and choose $T>0$ large enough so that, for any FSP we have $\|x(T)-x^*\| \le \epsilon$.
We claim that, as $r\to\infty$, 
\beql{eq-claim1}
\pr\{\|x^r(T)-x^*\|>2\epsilon\} \to 0,
\end{equation}
where the convergence is uniform across all initial states $x^r(0)$. This claim is true, because for an arbitrary sequence
of fixed initial states $x^r(0)$, we must have
$$
\limsup_{r\to\infty} \|x^r(T)-x^*\| \le \epsilon, ~~~\mbox{w.p.1},
$$
because w.p.1 we can always choose a subsequence of $r$ along which the 
u.o.c convergence to an FSP holds.
Claim \eqn{eq-claim1} implies the result,
because $\epsilon$ can be arbitrarily small.
$\Box$

\section{Original system. Optimality of a modified version of Greedy-D
%Proof of Theorem~\ref{th-fluid-stat-gt1}
}
\label{sec-proof-modified-gt1}

We now return to our original ``open'' system.
Unfortunately, the proof of Greedy-D optimality in the closed system does not 
carry over to the open system. The key reason can be informally described as follows.
If $v(t)$ is the exogenous arrival rate allocation vector in the fluid limit,
then the result analogous to Lemma~\ref{lem-F-decrease} no longer holds;
namely, property \eqn{eq-H-order} in its proof is not valid.
Suppose we have an edge $(\ell,j)$ on which the unique minimum 
$\min_{\ell'} \Delta_{\ell' j}$ is attained, $x_{\ell-e_j}=0$, $x_{\ell}>0$ and $\ell_j>0$.
There is the non-zero rate $\mu_j \ell_j x_{\ell}$ of departures along $(\ell,j)$.
In this case it is possible that $v_{\ell j} < \mu_j \ell_j x_{\ell}$, because
some type $j$ exogenous arrivals will find no servers in configuration
$\ell-e_j$. (In the closed system we must have
$v_{\ell j} \ge \mu_j \ell_j x_{\ell}$, because all departures along $(\ell,j)$
have the option of ``coming back'' along the same edge.)
Therefore, the argument leading to \eqn{eq-H-order}, and in fact the property itself,
does not hold.

In this section we will prove 
the asymptotic optimality 
% Theorem~\ref{th-fluid-stat-gt1} 
of a slightly modified version of the 
Greedy-D algorithm, called Greedy-DM, 
which, in a sense, ``emulates'' the behavior of Greedy-D in
the corresponding closed system. Informally the key idea of Greedy-DM
is to make decisions 
about placements of new exogenous type $i$ arrivals a little ``in advance'', 
at the instants of type $i$ departures. This is achieved by using ``placeholders,''
called tokens: when a type $i$ departure occurs, we ``pretend'' that 
immediately a new type $i$ arrival occurs, decide which server this arrival would be placed into
according to the Greed-D rule, and place a type $i$ token in that server.
When new actual type $i$ customers arrive, they first seek and replace
type $i$ tokens if there are any; if no type $i$ token is available,  the customer is placed according to the Greedy-D.
(In addition to being replaced by actual customers, we make the tokens "impatient" -- they leave the system by themselves
after a random exponentially distributed time.) The analysis in this section shows that, in the fluid limit,
``all'' tokens are replaced by actual arrivals and ``all'' actual arrivals
replace tokens. This means that the allocation of actual customer arrival rates
is ``equal'' to that of tokens; the latter, in turn, satisfies same properties
as the rate allocation in the closed system.

\begin{definition}[Greedy-DM discipline] Suppose the weights 
$W_k(x)=x_k^{\alpha}$ are given, as in the
definition of standard Greedy-D. At any given time there are two kinds of type $i$ customers
-- actual customers being served as usual and {\em tokens}. 
%(Tokens, as shown below,
%serve as temporary ``placeholders'' for the actual customers arriving later.)
For the purposes of defining server configurations $k$ and the system state $X$,
type $i$ tokens are treated the same way as actual type $i$ customers.\\
When a departure of actual type $i$ customer from 
configuration $k$ at time $t$ occurs, 
the following actions are taken: a new token of type $i$
is immediately created; 
this new token is treated the same way as a new type $i$ arrival, and is placed for 
``service'' immediately; the token  is 
added to an available configuration $k$ 
(with either $k=0$ or $X_k(t-)>0$) such that $k+e_i\in \ck$ and the difference
$W_{k+e_i}(X(t-))-I\{k\ne 0\}W_k(X(t-))$ is the smallest (where $t-$ refers to the
time after the service completion, but before the token placement).\\
When a new exogenous arrival of an (actual) type $i$ customer occurs, it replaces
an arbitrarily chosen type $i$ token (anywhere in the system), 
if such is available; otherwise,
this arrival is 
added (as in the usual Greedy-D)
to an available configuration $k$ (with either $k=0$ or $X_k(t-)>0$) such that 
$k+e_i\in \ck$ and the difference
$W_{k+e_i}(X(t-))-I\{k\ne 0\}W_k(X(t-))$ is the smallest (where $t-$ refers to the
time just before the arrival).\\
Any token of any type anywhere in the system "completes service" at the rate $\mu_0>0$, independently of
anything else (i.e. the probability of service completion in a $dt$-long interval is $\mu_0 dt + o(dt)$);
``service completion'' of any token is treated 
the same way as service completion of an actual customer, except no new token is
created.
\iffalse
 each ``old'' token
of type $i$, anywhere in the system, independently ``completes service'' with the 
probability $1/(\hat Y_i(t-)+1)$, where $\hat Y_i(t-)$ is the number of {\em actual}
type $i$ customers in the system; 
%COM Maybe just make tokens to complete service at 
%a constant rate $\mu_{token}$. Probably simpler ENDCOM
``service completion'' of any token is treated 
the same way as service completion of an actual customer, except no new tokens are 
created and all other tokens are unaffected; in the special case when $\hat Y_i(t-)=0$,
the actions are as described above, except
no new type $i$ token is created or placed (and thus a departure that leaves no
actual type $i$ customers in the system also leaves no type $i$ tokens).\\
\fi
\end{definition}

The random process, describing the evolution of this system is more complex, but it is 
still an irreducible Markov chain. A {\em complete server configuration} is by definition
a pair $(k,\hat k)$, where vector $k=(k_1,\ldots,k_I)\in \ck$ gives the numbers all customers
(both actual and tokens) in a server, while vector $\hat k \le k$, $k\in \bar \ck$,
gives the number
of actual customers only. Therefore, the Markov process state at time $t$ is
the vector $\{X_{(k,\hat k)}(t)\}$, where the index $(k,\hat k)$ takes values
that are all possible complete server configurations, as described 
above.
Obviously, $\hat Y_i^r(t)\le Y_i^r(t)$ for all $i$ and $t$,
where $Y_i^r(t)$ and $\hat Y_i^r(t)$ is the total number of all and actual type $i$ customers,
respectively, and superscript $r$, as usual, indicates the system with parameter $r$.
Moreover, the behaviors of the processes $(Y_i^r(t), \hat Y_i^r(t)), ~t\ge 0$, are independent
across all $i$, with $\hat Y_i^r(\infty)$ 
having Poisson
distribution with mean $\rho_i r$.
Finally, by the Greedy-DM definition, at any time any existing type $i$ token
``completes service'' at the rate $\mu_0$.
Using these facts, we easily establish the following
\begin{lem}
\label{lem-complete-state-tight}
Markov chain 
$\{X_{(k,\hat k)}^r(t)\}, ~t\ge 0$,
is positive recurrent for each $r$. Moreover, the distributions of
$\{(\hat y_i^r(\infty), y_i^r(\infty)), ~i\in \ci\}= 
(1/r)\{(\hat Y_i^r(\infty), Y_i^r(\infty)), ~i\in \ci\}$
are tight, and
any limit in distribution \\
$\{(\hat y_i(\infty), y_i(\infty)), ~i\in \ci\}$ 
is such that $\hat y_i(\infty)=\rho_i$ and $y_i(\infty)\le \rho_i + \lambda_i /\mu_0$
for all $i$. Consequently, the distributions of 
$\{x_{(k,\hat k)}^r(\infty)\}= (1/r)\{X_{(k,\hat k)}^r(\infty)\}$
are tight.
\end{lem}

{\em Proof.} 
Consider the process $(Y_i^r(\cdot), \hat Y_i^r(\cdot))$ for a fixed $r$ and $i$.
Consider a different process $(\bar Y_i^r(\cdot), \hat Y_i^r(\cdot))$, defined as follows:
actual type $i$ customers arrive the same way, and depart after their service completion,
so that $\hat Y_i^r(t)$ is same as before;
upon every service completion of an actual $i$-customer, one type $i$ token is created,
which then stays in the system for a random time, exponentially distributed with 
mean $1/\mu_0$, independently of anything else, and then leaves the system;
$\bar Y_i^r(t)$ is the total number of all type $i$ customers (actual and tokens) in the
system. Obviously, processes $(Y_i^r(\cdot), \hat Y_i^r(\cdot))$ and 
$(\bar Y_i^r(\cdot), \hat Y_i^r(\cdot))$ can be constructed on a common probability
space so that $Y_i^r(t) \le \bar Y_i^r(t)$. But, $\bar Y_i^r(t)$ is simply the number
of customers in the infinite-server system with input rate $\lambda_i r$ and mean service 
time $1/\mu_i + 1/\mu_0$. Therefore, $(\bar Y_i^r(\cdot), \hat Y_i^r(\cdot))$
is positive recurrent; moreover, in stationary regime,
 $\bar Y_i^r(\infty)$ and $\hat Y_i^r(\infty)$
have Poisson distributions
with means $\lambda_i r (1/\mu_i + 1/\mu_0)$ and $\rho_i r$, respectively.
The lemma results easily follow.
$\Box$

\iffalse
is such that $y_i(\infty)=\hat y_i(\infty)=\rho_i$ for all $i$.
This immediately implies the following further properties.
The distributions of $[x_{(k,\hat k)}^r(\infty)]= (1/r)[X_{(k,\hat k)}^r(\infty)]$
are tight, and
any limit in distribution
$[x_{(k,\hat k)}(\infty)]$ is such that $x_{(k,\hat k)}(\infty)=0$ unless $\hat k=k$.
(Otherwise, $y_i(\infty)=\hat y_i(\infty)=\rho_i$ would be violated.
So, in the limit, the ``amount'' of configurations $(k,\hat k)$ with non-zero
number of tokens is zero.)
In particular, the distribution of $[x_{(k,\hat k)}(\infty)]$ is determined by
that of its projection $x(\infty)=[x_k(\infty), ~k\in\ck]$, 
where $x_k(\infty)=\sum_{\hat k \le k} x_{(k,\hat k)}(\infty)$;
moreover, $x(\infty)\in \cx$ w.p.1.
\fi

Let $\{x_{(k,\hat k)}\}$ denote a vector with non-negative components, with indices $(k,\hat k)$
being all possible complete configurations. Denote by $x=\{x_k\}$, $y=\{y_i\}$ and 
$\hat y=\{\hat y_i\}$ its projections, with components being
$$
x_k=\sum_{\hat k ~:~\hat k \le k} x_{(k,\hat k)}, ~~ y_i = \sum_{(k,\hat k)} k_i x_{(k,\hat k)} =
\sum_{k} k_i x_{k}, ~~ \hat y_i = \sum_{(k,\hat k)} \hat k_i x_{(k,\hat k)}.
$$
Denote by $\bar \cx$ the set of those values of $\{x_{(k,\hat k)}\}$ satisfying
condition
$$
y_i=\hat y_i=\rho_i, ~~\forall i.
$$
Obviously, for any $\{x_{(k,\hat k)}\} \in \bar \cx$, its $x$-projection is 
an element of $\cx$. Also, 
note that the condition $y_i=\hat y_i, ~\forall i,$ is equivalent to 
\beql{eq-diag}
x_{(k,\hat k)}=0 ~~\mbox{unless}~~ \hat k = k,
\end{equation}
and therefore \eqn{eq-diag} holds
 for any $\{x_{(k,\hat k)}\} \in \bar \cx$. 
\iffalse
Finally, denote by
$x^{**}$ the element of $\bar \cx$ such that its $x$-projection is 
the optimal point $x^*$ is defined in \eqn{eq-xstar}.
\fi

The main result of this section is

\begin{thm}
\label{th-fluid-stat-gt1-Greedy-DM}
Consider a sequence of original systems, indexed by $r$,
under the Greedy-DM algorithm with $\alpha>0$.
% with  $W_k(X)=X_k^\alpha$, $\alpha>0$.
Let $\{x_{(k,\hat k)}^r(\infty)\}$ denote the random (complete) 
state of the fluid-scaled process
in a stationary regime.
Then, as $r\to\infty$,
\beql{eq-close-to-diag}
d(\{x_{(k,\hat k)}^r(\infty)\},\bar\cx) \implies 0,
\end{equation}
\beql{eq-conv-to-optimal}
x^r(\infty) \implies x^*,
\end{equation}
%$$
%\{x_{(k,\hat k)}^r(\infty)\} \implies x^{**}.
%$$
where $x^*$ is defined in \eqn{eq-xstar}.
\end{thm}

For each $(k,i)\in \cm$ and $\hat k \le k$, consider independent
unit-rate Poisson processes $\hat \Pi_{(k,\hat k), i}(t), ~t\ge 0$,
and $\tilde \Pi_{(k,\hat k), i}(t), ~t\ge 0$; and for 
each $i\in \ci$ -- an independent unit-rate Poisson process $\hat \Pi_i(t), ~t\ge 0$.
Each of these processes satisfies the functional strong law of large numbers (FSLLN) 
analogous to \eqn{eq-flln-poisson}.

The Markov process $\{X^r_{(k,\hat k)}(\cdot)\}$ for each $r$ is driven by the common
set of independent
Poisson processes $\hat \Pi_{(k,\hat k), i}(\cdot)$, $\tilde \Pi_{(k,\hat k), i}(\cdot)$ and $\hat \Pi_{i}(\cdot)$, in the natural way,
as follows:
$$
\hat D^r_{(k,\hat k), i}(t) = 
\hat \Pi_{(k,\hat k), i} (\int_0^t X_{(k,\hat k)}^r(\xi) \hat k_i \mu_i d\xi),
$$
$$
\tilde D^r_{(k,\hat k), i}(t) = 
\tilde \Pi_{(k,\hat k), i} (\int_0^t X_{(k,\hat k)}^r(\xi) (k_i-\hat k_i) \mu_0 d\xi),
$$
$$
\hat A^r_i(t) = \hat \Pi_i(r \lambda_i t),
$$
where $\hat D^r_{(k,\hat k), i}(t)$ and $\tilde D^r_{(k,\hat k), i}(t)$
is the number of the type-i departures 
from the configuration $(k,\hat k)$, due to the service completions of
actual customers and tokens, respectively, and $\hat A^r_i(t)$ 
is the number of the exogenous type-i arrivals of actual customers;
all in the interval $[0,t]$. Clearly, the entire process sample path is
(w.p.1) uniquely determined by the initial state $\{X^r_{(k,\hat k)}(0)\}$,
the realizations of the driving Poisson processes
and the Greedy-DM discipline. In particular, the realizations of the following 
processes are uniquely determined:\\
the number of type-i departures from configuration $k$ due to actual and token
service completions, and their total
$$
\hat D^r_{k,i}(t) = \sum_{\hat k} \hat D^r_{(k,\hat k), i}(t), ~
\tilde D^r_{k,i}(t) = \sum_{\hat k} \tilde D^r_{(k,\hat k), i}(t), ~
D^r_{k,i}(t) = \hat D^r_{k,i}(t) + \tilde D^r_{k,i}(t);
$$
the number of type-i token ``arrivals'' allocated (upon type-i actual departures) along
edge $(k,i)$, $\tilde A^r_{k,i}(t)$; 
the number of type-i actual exogenous arrivals allocated along
edge $(k,i)$, {\em without replacing an existing i-token}, $\hat A^{**,r}_{k,i}(t)$
(such arrivals change the complete configuration from $(k-e_i,\hat k-e_i)$
to $(k,\hat k)$; the number of type-i actual exogenous arrivals,
{\em that replace an existing token in configuration $k$},
 $\hat A^{*,r}_{k,i}(t)$
(such arrivals change the complete configuration from $(k,\hat k-e_i)$
to $(k,\hat k)$ -- so these do {\em not} count as arrivals into $k$); 
the total number of all type-i arrivals into configuration $k$,
$A^r_{k,i}(t)(t) = \tilde A^r_{k,i}(t) + \hat A^{**,r}_{k,i}(t)$.
The following relations obviously hold:
$$
\hat A^r_i(t) = \sum_{k:(k,i)\in\cm} [\hat A^{*,r}_{k,i}(t)+\hat A^{**,r}_{k,i}(t)],~~\forall i,
$$
\beql{eq-token-conserv}
\sum_{k:(k,i)\in\cm} \tilde A^r_{k,i}(t) = \sum_{k:(k,i)\in\cm} \hat D^r_{k,i}(t), ~~\forall i.
\end{equation}

We introduce fluid-scaled processes:
$$
x^r_{k}(t) = \frac{1}{r} X^r_{k}(t),~~x^r_{(k,\hat k)}(t) = \frac{1}{r} X^r_{(k,\hat k)}(t),
~~\hat a^r_{i}(t) = \frac{1}{r} \hat A^r_{i}(t),
$$
and similarly defined
$$
\hat d^r_{ki}(t), ~\tilde d^r_{ki}(t), ~ d^r_{ki}(t), ~ 
~\tilde a^r_{ki}(t), ~\hat a^{*,r}_{ki}(t), ~\hat a^{**,r}_{ki}(t), ~a^r_{ki}(t).
$$
%$$
%\hat d^r_{ki}(t) = \frac{1}{r} \hat D^r_{ki}(t), ~~a^r_{ki}(t) = \frac{1}{r} A^r_{ki}(t),
%$$

A set of Lipschitz continuous functions
$[\{x_k(\cdot)\}, \{x_{(k,\hat k)}(\cdot)\}, \{\hat a_i(\cdot)\}, \{\hat d_{ki}(\cdot)\},\ldots]$
on the time interval $[0,\infty)$ we call a {\em fluid sample path} (FSP), if there exist
realizations of driving Poisson processes, satisfying the FSLLN analogous to
\eqn{eq-flln-poisson}, 
and a fixed subsequence of $r$, along which
%\beql{eq-fsp-def-closed}
\begin{eqnarray}
& [\{x^r_k(\cdot)\}, \{x^r_{(k,\hat k)}(\cdot)\}, \{\hat a^r_i(\cdot)\}, \{\hat d^r_{ki}(\cdot)\},\ldots]
\to \nonumber \\
& [\{x_k(\cdot)\}, \{x_{(k,\hat k)}(\cdot)\}, \{\hat a_i(\cdot)\}, \{\hat d_{ki}(\cdot)\},\ldots],
~~u.o.c. \label{eq-fsp-def-Greedy-DM}
\end{eqnarray}
%\end{equation}

\begin{lem}
\label{lem-conv-to-fsp-Greedy-DM}
Suppose the initial states $\{x_{(k,\hat k)}^r(0)\}$ are fixed and are such that 
$\{x_{(k,\hat k)}^r(0)\}\to \{x_{(k,\hat k)}(0)\}$.
Then, w.p.1 for any subsequence of $r$ there exists a further
subsequence of $r$, along which the convergence \eqn{eq-fsp-def-Greedy-DM} holds,
where the limit is an FSP.
\end{lem}

{\em Proof} is, again, very standard -- it is a more general version of that
of Lemma~\ref{lem-conv-to-fsp-closed}. We omit details.
$\Box$

\iffalse
Essentially, it suffices to observe that, since each $x^r_k(t)\le 1$,
all functions $\int_0^t x^r_k(\xi)d\xi$, for all $r$, are uniformly Lipschitz,
and therefore there always exits a subsequence of $r$ along which they converge u.o.c.
to Lipschitz functions. Using this, along with \eqn{eq-fsp-def-Greedy-DM}, we
easily verify the convergence of all $x^r_k(\cdot),d^r_{ki}(\cdot),a^r_{ki}(\cdot)$,
along possibly a further subsequence,
and the fact that the limits are Lipschitz.
\fi

\begin{lem}
\label{lem-fsp-properties--Greedy-DM-y}
Consider an FSP with the initial state $\{x_{(k,\hat k)}(0)\}$.
Recall notation:
$$
y_i(t) = \sum_{(k,\hat k)} k_i x_{(k,\hat k)}(t), 
~~ \hat y_i(t) = \sum_{(k,\hat k)} \hat k_i x_{(k,\hat k)}(t),
$$
and denote $\tilde y_i(t) = y_i(t)- \hat y_i(t)$. 
Then, at any regular point $t$, for any $i$,
\beql{eq-yhat}
(d/dt) \hat y_i(t) = \lambda_i - \mu_i \hat y_i(t),
\end{equation}
\beql{eq-ytilde}
(d/dt) \tilde y_i(t) = \left\{ \begin{array}{ll}
               -\lambda_i + \mu_i \hat y_i(t) - \mu_0 \tilde y_i(t), 
                & \mbox{if}~\tilde y_i(t)>0\\
               \max\{0,-\lambda_i + \mu_i \hat y_i(t) - \mu_0 \tilde y_i(t)\}, 
                & \mbox{if}~\tilde y_i(t)=0
                \end{array}
       \right.
\end{equation}
In particular, for any $i$, the convergence
\beql{eq-y-converge}
(\hat y_i(t), \tilde y_i(t)) \to (\rho_i, 0), ~~\forall i,
\end{equation}
holds and is uniform in initial states $\{x_{(k,\hat k)}(0)\}$ from a compact set,
and
\beql{eq-y-equilibrium}
(\hat y_i(0), \tilde y_i(0))= (\rho_i, 0)~~\mbox{implies}~~ 
(\hat y_i(t), \tilde y_i(t))= (\rho_i, 0), ~\forall t.
\end{equation}
\end{lem}

{\em Proof.} Equation \eqn{eq-yhat} is very standard, describing an FSP for 
an $M/M/\infty$ system. Equation \eqn{eq-ytilde} for the $\tilde y_i(t)>0$ case
is also a very basic; \eqn{eq-ytilde} for the $\tilde y_i(t)=0$ case is easily
verified by considering the behavior of pre-limit trajectories is a small
interval $[t,t+\Delta t]$, and considering the three cases,
$-\lambda_i + \mu_i \hat y_i(t) - \mu_0 \tilde y_i(t)<0$, $=0$ and $>0$,
separately. (See e.g. \cite{St2001} for this type of argument in more detail.) We omit details.
$\Box$

\begin{lem}
\label{lem-fluid-stat-gt1-Greedy-DM}
Consider a sequence of original systems, indexed by $r$,
under the Greedy-DM algorithm.
% with  $W_k(X)=X_k^\alpha$, $\alpha>0$.
Let $\{x_{(k,\hat k)}^r(\infty)\}$ denote the random complete 
state of the fluid-scaled process
in a stationary regime.
% under the Greedy-D discipline
%with  $W_k(X)=X_k^\alpha$, $\alpha>0$.
Then, any subsequence of $r$ has a further subsequence, such that
$$
\{x_{(k,\hat k)}^r(\infty)\} \implies\{x_{(k,\hat k)}(\infty)\},
$$
where $\{x_{(k,\hat k)}(\infty)\}\in \bar \cx$ w.p.1.
%the limit is such that $x_{(k,\hat k)}(\infty)=0$ unless $\hat k=k$,
%and $x(\infty)=x^*$. ($x^*$ is defined in \eqn{eq-xstar}.)
\end{lem}

{\em Proof.} Fix arbitrary $\delta>0$, and a sufficiently large compact
$B$ so that for all large $r$, 
\beql{eq-temp-123}
\pr[\{x_{(k,\hat k)}^r(\infty)\}\in B ] \ge 1-\delta.
\end{equation}
(We can do that by Lemma~\ref{lem-complete-state-tight}.)
Fix arbitrary $\epsilon>0$ and choose $T>0$ large enough so that for any FSP 
with $\{x_{(k, \hat k)}(0)\}\in B$,
we have $d(\{x_{(k, \hat k)}(T)\},\bar\cx)\le \epsilon$. 
(We can do that by Lemma~\ref{lem-fsp-properties--Greedy-DM-y}.)
Fix arbitrary $\delta_1>0$. 
We claim that for all sufficiently large $r$,
\beql{eq-claim555}
\{x^r_{(k, \hat k)}(0)\} \in B ~~\mbox{implies}~~
\pr\{d(\{x^r_{(k, \hat k)}(T)\},\bar\cx)>2\epsilon \} \le \delta_1.
\end{equation}
This claim is true, because 
for an arbitrary sequence
of fixed initial states $\{x^r_{(k, \hat k)}(0)\} \in B$, we must have
$$
\limsup_{r\to\infty}  d(\{x_{(k, \hat k)}(T)\},\bar\cx)\le \epsilon, ~~~\mbox{w.p.1}.
$$
(This follows from Lemma~\ref{lem-conv-to-fsp-Greedy-DM}.)
By \eqn{eq-temp-123} and \eqn{eq-claim555}, 
for all large $r$, a stationary version of the process
is such that
$$
\pr\{d(\{x^r_{(k, \hat k)}(T)\},\bar\cx) \le 2\epsilon \} \ge (1-\delta)(1-\delta_2).
$$
Therefore, any limit-in-distribution $\{x_{(k, \hat k)}(T)\}$ is such that
$$
\pr\{d(\{x_{(k, \hat k)}(T)\},\bar\cx) \le 2\epsilon \} \ge 
\limsup \pr\{d(\{x^r_{(k, \hat k)}(T)\},\bar\cx) \le 2\epsilon \} \ge (1-\delta)(1-\delta_2).
$$
Since $\delta, \delta_2, \epsilon$ are arbitrary positive, 
$\pr[\{x_{(k, \hat k)}(T)\}\in \bar \cx] =1$.
$\Box$

\begin{lem}
\label{lem-fsp-properties--Greedy-DM}
Consider an FSP with the initial state $(x_{(k,\hat k)}(0))\in \bar \cx$.
% such that $y_i(0)=\hat y_i(0) = \rho_i$ for each $i$.
(In particular, $x(0)\in\cx$.)
Then $(x_{(k,\hat k)}(t))\in \bar \cx$ for all $t\ge 0$. In addition,
\iffalse
Then the following properties hold:
\beql{eq-y-equilibrium222}
y_i(t)=\hat y_i(t) \equiv \rho_i,
\end{equation}
in particular, 
\beql{eq-y-equilibrium333}
x_{(k,\hat k)}(t)=0 ~\mbox{unless}~ \hat k=k;
\end{equation}
\fi
at any regular point $t$,
using notation $w_{ki}(t)=(d/dt)d_{ki}(t)$, $\hat w_{ki}(t)=(d/dt)\hat d_{ki}(t)$,
$v_{ki}(t)=(d/dt)a_{ki}(t)$, $\tilde v_{ki}(t)=(d/dt)\tilde a_{ki}(t)$, we have
\beql{eq-w-what}
w_{ki}(t)=\hat w_{ki}(t)=k_i \mu_i x_k, ~~\forall (k,i)\in \cm,
\end{equation}
\beql{eq-what-vtilde}
\sum_{k:(k,i)\in \cm} \hat w_{ki}(t) = \sum_{k:(k,i)\in \cm} \tilde v_{ki}(t) = \lambda_i, ~~\forall i\in \ci,
\end{equation}
\beql{eq-v-vtilde}
v_{ki}(t) = \tilde v_{ki}(t), ~~\forall (k,i)\in \cm,
\end{equation}
\beql{eq-x-dynamics}
(d/dt) x_k(t) = [\sum_{i:k-e_i\in\bar\ck} v_{ki} - \sum_{i:k+e_i\in\ck} v_{k+e_i,i} ]
                     - [\sum_{i:k-e_i\in\bar\ck} w_{ki} - \sum_{i:k+e_i\in\ck} w_{k+e_i,i} ], ~~\forall k\in\ck,
\end{equation}
\beql{eq-F-decrease-Greedy-DM}
\frac{d}{dt}F(x(t)) = D(v(t),x(t)) \le D_{min}(x(t)).
\end{equation}
\end{lem}

{\em Proof.} Since we have \eqn{eq-y-equilibrium}, $\{x_{(k,\hat k)}(t)\}\in \bar \cx$ 
holds by definition.
%\eqn{eq-y-equilibrium222} is just a form of \eqn{eq-y-equilibrium}.
Relation \eqn{eq-w-what} holds because
$$
\hat w_{ki}(t)= \sum_{\hat k \le k} \hat k_i \mu_i x_{(k,\hat k)}(t), ~~
 w_{ki}(t)-\hat w_{ki}(t) = \sum_{\hat k \le k} (k_i-\hat k_i) \mu_0 x_{(k,\hat k)}(t),
$$
and $\{x_{(k,\hat k)}(t)\} \in \bar \cx$.
%\eqn{eq-y-equilibrium333}. 
We obtain \eqn{eq-what-vtilde} from
the limit form of \eqn{eq-token-conserv}, namely
$$
\sum_{k:(k,i)\in\cm} \tilde a_{k,i}(t) = \sum_{k:(k,i)\in\cm} \hat d_{k,i}(t),
$$
and from $\sum_k w_{ki}(t)=\sum_k \hat w_{ki}(t)= \mu_i \hat y_i(t)=\lambda_i$.
Relation \eqn{eq-v-vtilde} follows from the fact that $v_{ki}(t) \ge \tilde v_{ki}(t)$,
and the strict inequality cannot hold for any $(k,i)$, because otherwise we would have
for at least one $i$
$$
(d/dt) y_i(t) = \sum_k v_{ki}(t) - \sum_k w_{ki}(t) > 
\sum_k \tilde v_{ki}(t) - \sum_k \hat w_{ki}(t) = 0.
$$
Equation \eqn{eq-x-dynamics} is automatic: the RHS is just the difference between arrival
and departure rates to/from configuration $k$. Finally, since $v(t)=\tilde v(t)$,
and the rates $\tilde v_{ki}(t)$ are those of ``arriving'' tokens (which immediately
follow service completions of actual customers), the argument in the proof
of Lemma~\ref{eq-F-decrease-closed} (for the closed system) applies,
and we obtain \eqn{eq-F-decrease-Greedy-DM}.
$\Box$

\iffalse
The rate at which new type $i$ tokens are created 
(in the FSP under consideration) 
is equal to the actual customer departure rate $\rho_i \mu_i=\lambda_i$.
The rate at which type $i$ tokens ``complete service'' is zero, simply
because their amount $y_i(t) - \hat y_i(t)$ is zero. Therefore, the rate at which type $i$ 
tokens are replaced by new arrivals of the actual type $i$ customers is $\lambda_i$.
But this exactly matches the rate of new actual arrivals. We conclude that,
in the FSP, new type $i$ actual arrivals allocation $v(t)$ is equal to the allocation
of type $i$ tokens. We repeat the argument in the proof of Lemma~\ref{lem-F-decrease},
to complete this proof. 
$\Box$
\fi

As a corollary we obtain the analog of Lemma~\ref{lem-fluid-convergence-closed}.

\begin{lem}
\label{lem-fluid-convergence-DM}
Consider an FSP (for the original system under Greedy-DM algorithm)
with the initial state $\{x_{(k,\hat k)}(0)\} \in \bar \cx$.
% such that
%$y_i(0)=\hat y_i(0) = \rho_i$ for each $i$.
(In particular, $x(0)\in\cx$.) 
Then
\beql{eq-fluid-convergence-DM}
x(t)\to x^*.
\end{equation}
The convergence is uniform across all initial states in $\bar \cx$.
\end{lem}

{\em Proof of Theorem~\ref{th-fluid-stat-gt1-Greedy-DM}.} 
Convergence \eqn{eq-close-to-diag} has already been proved in 
Lemma~\ref{lem-fluid-stat-gt1-Greedy-DM}.
We fix $\epsilon>0$ and choose $T>0$ large enough so that for any FSP 
with $\{x_{(k, \hat k)}(0)\}\in\bar\cx$,
we have $\|x(T)-x^*\| \le \epsilon$. 
%$[x^*_{(k, \hat k)}]$ the vector such that $x^*_{(k, k)}=x^*_k$ and $x^*_{(k, k)}=0$
%when $\hat k \ne k$.
We claim that for any $\delta_1>0$ there exists a sufficiently small  $\delta_2>0$ 
such that for all sufficiently large $r$,
\beql{eq-claim111}
d(\{x^r_{(k, \hat k)}(0)\},\bar\cx)\le \delta_2 ~~\mbox{implies}~~
\pr\{\|x^r(T)-x^*\|>2\epsilon \} \le \delta_1,
\end{equation}
This claim is true, because 
for an arbitrary sequence
of fixed initial states $\{x^r_{(k, \hat k)}(0)\}\to \bar\cx$, we must have
$$
\limsup_{r\to\infty} \|x^r(T)-x^*\| \le \epsilon, ~~~\mbox{w.p.1}.
$$
Constants $\epsilon$ and $\delta_1$ can be arbitrarily small;
we also know that for any $\delta_2$, 
$\pr\{d(\{x^r_{(k, \hat k)}(\infty)\},\bar\cx)\le \delta_2\}\to 1$ as $r\to\infty$.
Therefore, claim \eqn{eq-claim111} implies \eqn{eq-conv-to-optimal}.
$\Box$

\section{The case of vector-packing constraints \eqn{eq-intro1}: 
Greedy algorithm with aggregate configurations}
\label{sec-aggregation}

So far in the paper we did not exploit a possible underlying structure of
packing constraints. Instead, we worked with a formally defined set $\bar \ck$ 
of possible configurations. Now we will consider a special case: suppose the 
configuration set $\bar \ck$ is defined by vector packing constraints 
\eqn{eq-intro1}. 

We say that two configurations $k$ and $k'$ are equivalent if they require 
same total amounts of resources of each type:
\beql{eq-intro1-222}
\sum_i k_i b_{i,n} = \sum_i k'_i b_{i,n}, ~~\forall n.
\end{equation}
A class of equivalent configurations is denoted $q$; we will call it {\em aggregate
configuration}, or {\em a-configuration}.
 Zero a-configuration, denoted (with some notation abuse) 
by $q=0$, is the one containing the sole configuration $k=0$; 
by convention $X_0=0$, where the subscript $0$ can refer to either zero configuration 
or zero a-configuration.
The sets of all a-configurations and non-zero a-configurations
are denoted by $\bar \cq$ and $\cq$, respectively. 
We write $q(k)$ for the aggregate configuration containing $k$.
We use notation
$$
X_q = \sum_{k\in q} X_k,
$$
and similarly for other quantities summed up over an a-configuration $q$. 
Clearly, vector $\{X_q, ~q\in \cq\}$ is a projection
of $X = \{X_k, ~k\in \ck\}$. 

In this section we show that (versions of) the Greedy algorithm, using quantities $X_q$ instead of $X_k$,
asymptotically minimizes
$$
\sum_{q\in\cq} X_q^{1+\alpha}(\infty),
$$
which, again, approximates (when $\alpha>0$ is small) the total number of
occupied servers $\sum_{q} X_q(\infty)\equiv \sum_k X_k(\infty)$ in a stationary regime.
The difference, however, is that $|\cq|$ can be much smaller than $|\ck|$,
thus making the Greedy algorithm easier to implement in practice.

\subsection{Results.} 

For $x\in\R_+^{|\ck|}$ consider function 
$$
\Phi(x)=\sum_q (1+\alpha)^{-1} x_q^{1+\alpha} \equiv 
\sum_q (1+\alpha)^{-1} [\sum_{k\in q} x_k]^{1+\alpha}
$$ 
with parameter $\alpha>0$,
and the following the convex optimization problem 
\beql{eq-opt-greater1-aaa}
\min_{x\in\R_+^{|\ck|}} \Phi(x)
\end{equation}
subject to
\beql{eq-cons-laws222-aaa}
\sum_{k\in\ck} k_i x_k = \rho_i, ~~\forall i.
\end{equation}
%\beql{eq-nonneg222}
%x_k \ge 0, ~~\forall k\in \ck,
%\end{equation}
We denote by $\cx^*$ the set of optimal solutions of this problem; obviously, 
$\cx^* \subseteq \cx$.

\begin{definition}[Greedy discipline with a-configurations]  $ $ 
\begin{enumerate}

\item {\em Integral form (Greedy-I-AC)}. 
A type $i$ customer arriving at time $t$ is 
added to an available a-configuration $q$ (with either $q=0$ or $X_q(t-)>0$)
such that the addition does not violate the vector packing constraints and
the increment $\Phi(X(t))-\Phi(X(t-))$ is the smallest.
The ties {\em between a-configurations}
are broken according to an arbitrary deterministic rule.
The choice of a server {\em within} the chosen a-configuration is random uniform.

\item {\em Differential form (Greedy-D-AC)}. For each $q\in\cq$,
denote
$W_q(x)=(\partial / \partial x_q) \Phi(x)=x_q^{\alpha}, ~x\in\R_+^{|\ck|}$. 
A type $i$ customer arriving at time $t$ is 
added to an available configuration $k$ (with either $q=0$ or $X_q(t-)>0$) 
such that the addition does not violate the vector packing constraints and the difference
$W_{q+e_i}(X(t-))-I\{q\ne 0\}W_q(X(t-))$ is the smallest. 
[Here $q+e_i$ denotes the a-configuration containing configurations $k+e_i, ~k\in q$
(and possibly other configurations).]
The ties {\em between a-configurations}
are broken according to an arbitrary deterministic rule.
The choice of a server {\em within} the chosen a-configuration is random uniform.

\end{enumerate}

\end{definition}

\begin{thm}
\label{th-fluid-stat-gt1-closed-aaa}
Consider a sequence of closed systems, indexed by $r$,
and let $x^r(\infty)$ denote the random state of the (fluid-scaled) process
is a stationary regime, under the Greedy-D-AC discipline with $\alpha>0$.
Then, as $r\to\infty$,
$$
d(x^r(\infty),\cx^*) \implies 0.
$$
\end{thm}

\begin{definition}[Greedy-DM-AC discipline] This discipline, which uses tokens,
is the modification of Greedy-D-AC, completely analogous to the modification
of Greedy-D that leads to Greedy-DM.
\end{definition}

\begin{thm}
\label{th-fluid-stat-gt1-Greedy-DM-aaa}
Consider a sequence of original systems, indexed by $r$,
under the Greedy-DM-AC algorithm with $\alpha>0$.
% with  $W_k(X)=X_k^\alpha$, $\alpha>0$.
Let $\{x_{(k,\hat k)}^r(\infty)\}$ denote the random (complete) 
state of the fluid-scaled process
in a stationary regime.
Then, as $r\to\infty$,
$$ %\beql{eq-close-to-diag}
d(\{x_{(k,\hat k)}^r(\infty)\},\bar\cx) \implies 0,
$$ %\end{equation}
$$ %\beql{eq-conv-to-optimal}
d(x^r(\infty),\cx^*) \implies 0.
$$ %\end{equation}
\end{thm}

In the rest of this section we will 
consider the closed system under Greedy-D-AC, and will
prove Theorem~\ref{th-fluid-stat-gt1-closed-aaa}.
We will omit the proof of Theorem~\ref{th-fluid-stat-gt1-Greedy-DM-aaa} 
which is ``obtained from'' that of Theorem~\ref{th-fluid-stat-gt1-closed-aaa}
in exactly same way as the proof of Theorem~\ref{th-fluid-stat-gt1-Greedy-DM}
 was obtained from that of Theorem~\ref{th-fluid-stat-gt1-closed}.

\subsection{Optimal set characterization and related properties.}

Using Lagrange multipliers $\eta_i$ for the constraints \eqn{eq-cons-laws222-aaa},
the Lagrangian of the problem \eqn{eq-opt-greater1-aaa}-\eqn{eq-cons-laws222-aaa} is
$$
\sum_q \frac{1}{1+\alpha} [\sum_{k\in q} x_k]^{1+\alpha} 
+ \sum_i \eta_i [\rho_i - \sum_k k_i x_k].
$$
We obtain the following characterization: vector $x\in \cx$
is an optimal solution of
\eqn{eq-opt-greater1-aaa}-\eqn{eq-cons-laws222-aaa}
(i.e. $x\in \cx^*$)
if and only if there exist constants
$\eta_i$ such that (using notation $u_k=\max\{\sum_i k_i \eta_i,0\}$)
\beql{x-star-characterize-aaa}
x_q^{\alpha} = \max_{k\in q} u_k,
%[\sum_{k\in q} x_k]^{\alpha} = \max_{k\in q} \max\{\sum_i k_i \eta_i,0\},
~~\forall q\in \cq,
\end{equation}
\beql{x-star-characterize-aaa3}
\left(u_k < \max_{k'\in q(k)} u_{k'} ~~\mbox{imples}~~ x_k=0\right), ~~\forall k\in \ck.
\end{equation}

More notation. Consider the following order relation on $\bar \cq$:
$q' \le q$ if $k' \le k$ for some $k'\in q'$ and $k\in q$.
$q'<q$ means $q'\le q$ and $q'\ne q$.
If a-configuration $q$ contains at least one $k$ with $k_i>0$,
we use a (slightly abusive) notation $q-e_i$ for the a-configuration
containing $k-e_i$, i.e. $q-e_i \doteq \{k - e_i ~|~ k\in q, k_i>0\}$;
otherwise, $q-e_i=\emptyset$. 
%COM need $q+e_i$?
Denote by $\cm^a$ the set of pairs $(q,i)$ such that
$q-e_i \ne \emptyset$.
For $x \in \R_+^{|\ck|}$ and $(q,i)\in \cm^a$ denote
$$
\Delta_{qi} = \Delta_{qi}(x) = x_q^{\alpha} - x_{q-e_i}^{\alpha}.
$$
%where, by convention, $x_k=0$ when $k=0$. 

\begin{lem}
\label{lem-NSI}
Consider the following property of an element $x\in \cx$.
(We will refer to it as {\em NSI-property} -- ``No Simple Improving allocation'').
For any 
two elements $(q,i), (q',i) \in \cm^a$ (with $q\ne q'$, but a common $i$),
condition 
\beql{eq-Delta-ineq-aaa}
\Delta_{q'i} < \Delta_{qi}
\end{equation}
implies either
\beql{eq-Delta-NSI-1}
x_k=0 ~~\mbox{for all $k\in q$ such that $k_i>0$}
\end{equation}
or
\beql{eq-Delta-NSI-2}
q'-e_i \ne 0 ~~\mbox{and}~~x_{q'-e_i}=0.
\end{equation}
If $x\in \cx$ satisfies the NSI-property, then condition \eqn{x-star-characterize-aaa}
holds.
\iffalse
(i) necessary for  $x\in \cx$ to be an optimal solution of
\eqn{eq-opt-greater1-aaa}-\eqn{eq-cons-laws222-aaa};\\
\fi
\end{lem}

{\em Proof.}
\iffalse
Statement (i) is obvious -- if NSI-property does not hold we can change $x$
so that the value of $\Phi(x)$ strictly decreases. 
Let us prove (ii). 
\fi
Consider $x\in \cx$ satisfying the NSI-property. 
For each $i$ define
$$
\underline \xi_i = \min_{k:~k_i>0, ~x_k>0} \Delta_{q(k),i}, ~~~
\overline \xi_i = \max_{k:~k_i>0, ~x_k>0} \Delta_{q(k),i}.
$$
It is easy to check that we cannot have $\overline \xi_i>0$
and $\underline \xi_i \le 0$, because this would violate the NSI-property.
Then, we can further define
\beql{eq-def-eta}
\eta_i= \left\{ \begin{array}{ll}
               \overline \xi_i, & \mbox{if}~\overline \xi_i>0\\
               \underline \xi_i, & \mbox{if}~\underline \xi_i \le 0
                \end{array}
       \right.
\end{equation}
Denote by $\ci^+$ the subset of those $i$ with $\eta_i>0$, 
and by $\ci^- = \ci \setminus \ci^+$ the remaining subset.
It is easy to check that for any fixed $i\in \ci^-$, 
we must have 
\beql{eq-equal555}
\Delta_{q(k),i}=\eta_i ~~\mbox{for all $k$ such that $k_i>0, ~x_k>0$},
\end{equation}
otherwise a contradiction to NSI-property is obtained.

Using notation $u_k=\max\{\sum_i k_i \eta_i,0\}$, let 
us define the values $x_q^0$ via
$$
[x_q^0]^{\alpha} = \max_{k\in q} u_k.
$$
It is easy to check that 
\beql{eq-ineq666}
[x_{q}^0]^{\alpha}-[x_{q-e_i}^0]^{\alpha} \ge \eta_i, ~~\forall (q,i)\in \cm^a.
\end{equation}

Let us prove \eqn{x-star-characterize-aaa}, which is equivalent to $x_q=x_q^0, ~\forall q$.
Suppose this is not true. Consider a minimal $q$ for which $x_q \ne x_q^0$.\\
Case (c1): suppose $x_q > x_q^0$. Then necessarily $x_q > 0$. Consider any $k\in q$
with $x_k>0$. \\
Sub-case (c1.1): suppose $k_i>0$ for some $i\in \ci^-$. Fix this $i$ and denote
$q'=q-e_i$. If $q'=0$, we obtain a contradiction to the definition of $\eta_i$,
and so $q' \ne 0$ must hold. Then $x_{q'}= x_{q'}^0$ by definition of $q$ (as a minimal
counterexample). 
Using \eqn{eq-ineq666}
%Since $[x_{q}^0]^{\alpha}-[x_{q'}^0]^{\alpha} \ge \eta_i$, 
we obtain
$x_{q}^{\alpha}-x_{q'}^{\alpha} = \Delta_{q,i}> \eta_i$ -- 
a contradiction to \eqn{eq-equal555}.
%(recall that since $i\in \ci^-$, $k_i>0, ~x_k>0$, 
%this difference must be equal to $\eta_i$). 
Thus sub-case (c1.1)
is impossible.\\
Sub-case (c1.2): suppose $k_i>0$ implies $i\in \ci^+$ and then $\eta_i>0$.
Fix an $i$ with $k_i>0$, and consider $q'=q-e_i$. If $q'=0$, 
we obtain a contradiction to the definition of $\eta_i$,
and so $q' \ne 0$ must hold. We have $x_{q'}= x_{q'}^0$ by definition of $q$;
%and $[x_{q}^0]^{\alpha}-[x_{q'}^0]^{\alpha} \ge \eta_i$; 
therefore, using \eqn{eq-ineq666},
$x_{q}^{\alpha}-x_{q'}^{\alpha} = \Delta_{q,i}> \eta_i$ -- a contradiction with
the definition of $\eta_i$. Sub-case (c1.2), and then case (c1), is impossible.\\
Case (c2): suppose $x_q < x_q^0$. Then necessarily $x_q^0 > 0$.
Let us fix a $k\in q$, on which the $\max_{k\in q} u_k >0$ is attained,
and so $u_k>0$. \\
Sub-case (c2.1): suppose $k_i>0$ for some $i\in \ci^-$. Fix this $i$ and denote
$q'=q-e_i$.
We cannot have $q'=0$, because that would imply $\eta_i>0$. Therefore, $q'\ne 0$
and $[x_{q}^0]^{\alpha}-[x_{q'}^0]^{\alpha} = \eta_i$, implying in particular
$x_{q'}^0>0$. 
But, $x_{q'}=x_{q'}^0$ and therefore $x_{q}^{\alpha}-x_{q'}^{\alpha} = \Delta_{q,i}< \eta_i$.
Recalling the definition of $\eta_i$, we obtain a contradiction to NSI-property.
Thus, sub-case (c2.1) is impossible.\\
Sub-case (c2.2): suppose $k_i>0$ implies $i\in \ci^+$ and then $\eta_i>0$.
Fix an $i$ with $k_i>0$, and consider $q'=q-e_i$. If $q'=0$, 
we have $\Delta_{q,i}< \eta_i$ -- a contradiction to NSI-property.
Therefore, $q'\ne 0$ and $x_{q'}=x_{q'}^0>0$ (because $\eta_i>0$ for all $i$ with $k_i>0$).
Then, $\Delta_{q,i}< \eta_i$ and we, again, obtain a contradiction to NSI-property.
Sub-case (c2.2), and then case (c2), is impossible.\\
The proof of \eqn{x-star-characterize-aaa} is complete.
$\Box$

\subsection{Fluid sample paths.}

We now define fluid sample paths for the closed system under Greedy-D-AC algorithm.
First, we will specify the construction of the process itself.
In addition to the set of unit-rate Poisson processes, driving
the service completions, we define primitive processes (common for each $r$),
driving the random uniform assignment of customers ``within'' each 
a-configuration $q$. Namely, for each $q$ we define an i.i.d. sequence
$\xi_q(1), \xi_q(2), \ldots$ of random variables, uniformly distributed in $[0,1]$.
The configurations $k\in q$ are indexed by $1,2,\ldots,|q|$ 
(in arbitrary fixed order). When an $m$-th customer of any type is
assigned to a-configuration $q$
(with $m$ referring to the order of assignment since initial time $0$,
and not to the customer type), this customer is assigned 
to a server in configuration $k'$ indexed by $1$ if 
$$
\xi_q(m) \in [0,X^r_{k'}/ X^r_{q}],
$$
it is assigned to
a server in configuration $k''$ indexed by $2$ if 
$$
\xi_q(m) \in (X^r_{k'}/ X^r_{q},(X^r_{k'}+X^r_{k''})/ X^r_{q}],
$$
and so on. (Note that necessarily $X^r_{q}>0$ -- otherwise there would be
no assignment to a-configuration $q$.) Denote
$$
g^r_q(s,\zeta) \doteq \sum_{m=1}^{\lfloor rs \rfloor} I\{\xi_q(m) \le \zeta\},
$$
where $s\ge 0$, $0\le \zeta \le 1$, and $\lfloor \cdot \rfloor$ denotes the integer part of a number. 
Obviously, from the strong law of large 
numbers (SLLN) and the monotonicity of $g^r_q(s,\zeta)$ on both arguments, we 
have the following functional SLLN
\beql{eq-flln-tie-break}
g^r_q(s,\zeta) \to s\zeta, ~~~\mbox{u.o.c.} ~~~\mbox{w.p.1}
\end{equation}
Clearly, the realization of the process is uniquely determined by the
initial state and the realizations of driving processes
$\Pi_{ki}(\cdot)$ and $(\xi_q(1), \xi_q(2), \ldots)$.

A set of Lipschitz continuous functions
$[\{x_k(\cdot),~k\in \ck\}, \{d_{ki}(\cdot),~(k,i)\in \cm\},\{a_{ki}(\cdot),~(k,i)\in \cm\}]$
on the time interval $[0,\infty)$ we call a {\em fluid sample path} (FSP), if there exist
realizations of $\Pi_{ki}(\cdot)$ satisfying \eqn{eq-flln-poisson},
realizations of $(\xi_q(1), \xi_q(2), \ldots)$ satisfying 
\eqn{eq-flln-tie-break},
and a fixed subsequence of $r$, along which convergence 
\eqn{eq-fsp-def-closed} holds. 

It is easy to see that the family of all FSPs is {\em uniformly} Lipschitz.

\iffalse
%\beql{eq-fsp-def-closed}
\begin{eqnarray}
& [\{x_k^r(\cdot),~k\in \ck\}, \{d_{ki}^r(\cdot),~(k,i)\in \cm\},\{a_{ki}^r(\cdot),~(k,i)\in \cm\}]
\to \nonumber \\
& [\{x_k(\cdot),~k\in \ck\}, \{d_{ki}(\cdot),~(k,i)\in \cm\},\{a_{ki}(\cdot),~(k,i)\in \cm\}],
~~u.o.c. \label{eq-fsp-def-closed}
\end{eqnarray}
\fi

We can easily verify that Lemmas~\ref{lem-conv-to-fsp-closed} and 
\ref{lem-fsp-properties-closed-basic} hold as is for Greedy-D-AC algorithm. 
Further, the following lemma
is analogous to Lemma~\ref{lem-F-decrease} (and has essentially same proof).

\begin{lem}
\label{lem-F-decrease-aaa}
Any FSP 
%(for the closed system under Greedy-D-AC) 
is such that
at any regular point $t$, 
%the actual allocation $v=v(t)$ is such that
\beql{eq-F-decrease-closed-aaa}
\frac{d}{dt}F(x(t)) \le 0.
\end{equation}
Moreover, unless $x(t)$ satisfies NSI-condition, the inequality 
\eqn{eq-F-decrease-closed-aaa} is strict.
\end{lem}

We will need the following FSP property, which follows from the
random uniform rule of Greedy-D-AC for assignments within each a-configuration.

\begin{lem}
\label{lem-zero-arrivals}
Consider an FSP. Suppose that at some $t> 0$, $x=x(t)$ is such that
for some $k\in \ck$ and $i$ we have: $x_k>0$, $k_i>0$,
$k'=k-e_i \ne 0$, $x_{k'}=0$, 
and $x_{q'}>0$ where $q'=q(k')$. Then $t$ is not a regular point.
\end{lem}

{\em Proof.} Suppose $t$ is a regular point.
We must have $\sum_{i'} v_{k'+e_{i'},i'}(t)=0$. Indeed, in a small interval $[t,t+\delta]$,
for all sufficiently large $r$,
the pre-limit sample paths defining the FSP are such that 
$$
\frac{x^r_{k'}}{x^r_{q'}} < \frac{C\delta}{x_{q'}-C\delta},
$$
where $C>0$ is some constant (depending on the Lipschitz constants for FSP components).
This, along with \eqn{eq-flln-tie-break}, implies that the fraction 
of the customers added in $[t,t+\delta]$
to servers in configuration $k'$, among those added to
a-configuration $q'$, is upper bounded by the RHS, which  
can be made arbitrarily small by choosing sufficiently small 
$\delta$. 
This means that $\sum_{i'}(a_{k'+e_{i'},i'}(t+\delta)-a_{k'+e_{i'},i'}(t))/\delta \downarrow 0$ as $\delta\to 0$,
which implies $\sum_{i'} v_{k'+e_{i'},i'}(t)=0$. Since $x_{k'}=0$,
obviously, $\sum_{i'} w_{k',i'}(t)=0$. However, $w_{k,i}(t)=\mu_i k_i x_k >0$.
Therefore, $(d/dt)x_{k'}(t) > 0$. This is a contradiction,
because if $t$ is regular, $x_{k'}(t) = 0$ implies $(d/dt)x_{k'}(t) = 0$.
$\Box$

We will also need the {\em continuity} and {\em shift} properties 
of the FSPs, which are quite generic. (See Sections 5 and 6 in \cite{St95}.
Although our model is different, essentially same proofs as in \cite{St95}
apply.) The {\em time shift by $\theta\ge 0$}, applied to an FSP
$[\{x_k(\cdot)\}, \{d_{ki}(\cdot)\},\{a_{ki}(\cdot)\}]$,
produces the set of functions with the same time argument $t\ge 0$, but
with $x_k(t)$ replaced by $x_k(\theta+t)$,
$d_{ki}(t)$ replaced by $d_{ki}(\theta+t)-d_{ki}(\theta)$,
$a_{ki}(t)$ replaced by $a_{ki}(\theta+t)-a_{ki}(\theta)$.

\begin{lem}
\label{lem-fsp-basic-aaa}
The family of FSPs satisfies the following properties.\\
(i) Continuity: If there is a converging sequence of FSPs,
indexed by $\beta$, namely
$$
[\{x_k^{(\beta)}(\cdot)\}, \{d_{ki}^{(\beta)}(\cdot)\},\{a_{ki}^{(\beta)}(\cdot)\}]
\to [\{x_k(\cdot)\}, \{d_{ki}(\cdot)\},\{a_{ki}(\cdot)\}], ~~\mbox{u.o.c.},
$$ 
then the limit is also an FSP.\\
(ii) Shift (or ``Memoryless''): The time shift of any FSP by any $\theta\ge 0$
is also an FSP.
\end{lem}

{\em Proof.}
(i) For each fixed index $\beta$, and the FSP associated with it,
consider a sequence of
(scaled) sample paths of the process, that define this FSP:
$$
[\{x_k^{(\beta,r)}(\cdot)\}, \{d_{ki}^{(\beta,r)}(\cdot)\},\{a_{ki}^{(\beta,r)}(\cdot)\}]
\to
[\{x_k^{(\beta)}(\cdot)\}, \{d_{ki}^{(\beta)}(\cdot)\},\{a_{ki}^{(\beta)}(\cdot)\}],
~~\mbox{u.o.c.}, ~~\mbox{as}~r\to\infty.
$$
Then, we can choose a subsequence of $r$, and the corresponding $\beta=\beta(r)$,
so that
$$
[\{x_k^{(\beta(r),r)}(\cdot)\}, \{d_{ki}^{(\beta(r),r)}(\cdot)\},\{a_{ki}^{(\beta(r),r)}(\cdot)\}]
\to
[\{x_k(\cdot)\}, \{d_{ki}(\cdot)\},\{a_{ki}(\cdot)\}],
~~\mbox{u.o.c.}, 
$$
and therefore the limit satisfies the definition of an FSP.\\
(ii) We pick a sequence of
(scaled) sample paths of the process, that define the FSP.
It is easy to see that the time shifts of these sample paths define the
FSP which is the time shift of the original one.
$\Box$

We are now in position to prove the following lemma, which is key (along with 
Lemmas~\ref{lem-NSI} and \ref{lem-zero-arrivals}) in our analysis of Greedy-D-AC algorithm.

\begin{lem}
\label{lem-drift-aaa}
Consider an FSP. Suppose $t$ is a regular point and $x(t)\not\in \cx^*$.
Then
\beql{eq-drift-aaa}
(d/dt) \Phi(x(t)) < 0.
\end{equation}
\end{lem}

{\em Proof.}
Suppose not, namely $(d/dt) \Phi(x(t)) = 0$. Then
for $x=x(t)$ the NSI-condition holds, and therefore 
\eqn{x-star-characterize-aaa} holds as well.
We will obtain a contradiction.
Condition \eqn{x-star-characterize-aaa3} does {\em not} hold (otherwise,
\eqn{x-star-characterize-aaa} and \eqn{x-star-characterize-aaa3}
would imply $x\in \cx^*$). Then, consider a minimal a-configuration $q$ for
which  \eqn{x-star-characterize-aaa3} is violated, namely:
for some $k\in q$,
\beql{x-star-characterize-counter}
u_k < \max_{k'\in q} u_{k'}, ~~ x_k>0.
\end{equation}
We will show that this is impossible. First, obviously, $x_q>0$.\\
Case (c3.1): suppose $k_i>0$ for some $i\in \ci^-$. Fix this $i$ and consider
$q'=q-e_i$. If $q'=0$, we have $\Delta_{q,i}>0$, which means 
$\eta_i$ cannot be negative -- a contradiction. Therefore,  $q'\ne 0$ and we 
must have $x_{q'}>0$ (because otherwise $\Delta_{q,i}>0$ leads, again, to 
the contradiction with $\eta_i<0$). Consider the set
$p\doteq \{k'+e_i ~|~ k'\in q'\}\subseteq q$; i.e., these are the configurations
in $q$ that are obtained by adding one type $i$ customer to configurations
in $q'$ -- it may or may not be a strict subset of $q$.
If $\max_{k'\in p} u_{k'} < \max_{k'\in q} u_{k'}$, then,
since $\max_{k'\in q'} u_{k'} + \eta_i \le \max_{k'\in p} u_{k'}$, we obtain
$\Delta_{q,i}>\eta_i$, which, along with $x_k>0$, leads to the contradiction 
with NSI-property. Therefore, $\max_{k'\in p} u_{k'} = \max_{k'\in q} u_{k'}$.
Then, there exists $k''\in \argmax_{k'\in q'} u_{k'}$ such that $x_{k''}>0$ and
$k''+e_i\in \argmax_{k'\in q} u_{k'}$. (Here we used the fact that for any 
$k'\in q'$, condition \eqn{x-star-characterize-aaa3} does hold --
recall that $q$ is a minimal a-configuration for which \eqn{x-star-characterize-aaa3}
is violated.) Note also that $u_{k-e_i} < \max_{k'\in q'} u_{k'} = u_{k''}$
and $x_{k-e_i} = 0$ (otherwise, again, \eqn{x-star-characterize-aaa3}
would be violated at $q'<q$).
We see that $x$ and the edge $(k,i)$ satisfy conditions
of Lemma~\ref{lem-zero-arrivals}. This means $t$ cannot be regular.
Thus, case (c3.1) is impossible.\\
Case (c3.2): suppose for any $i$ with $k_i>0$ we have $i\in \ci^+$. 
Fix one such $i$ and consider $q'=q-e_i$. 
%We will obtain a contradiction using the
%argument very similar to that in the case (c3.1).
If $q'=0$, we have $\Delta_{q,i}>\eta_i>0$  -- a contradiction with the definition of 
$\eta_i$. Therefore,  $q'\ne 0$ and then we must have $x_{q'}>0$ (because $\eta_j>0$
for each $j$ with $k_j>0$). Consider the set
$p\doteq \{k'+e_i ~|~ k'\in q'\}\subseteq q$. Note that $\max_{k'\in q'} u_{k'}>0$
and $\max_{k'\in q'} u_{k'}+\eta_i =\max_{k'\in p} u_{k'}>0$.
If $\max_{k'\in p} u_{k'} < \max_{k'\in q} u_{k'}$, then $\Delta_{q,i}>\eta_i>0$, 
which (along with $x_k>0$) contradicts the definition of $\eta_i$.
Therefore, $\max_{k'\in p} u_{k'} = \max_{k'\in q} u_{k'}$.
From this point on, the argument leading to a contradiction repeats that
in the case (c3.1) verbatim. Thus, the case (c3.2) is impossible.
The proof is complete.
$\Box$

\begin{lem}
\label{lem-uniform-drift-aaa}
For any $T>0$ and $\epsilon>0$, there exists $\delta>0$ such that the following property
holds uniformly on all FSPs and all $t_0\ge 0$:
\beql{eq-uniform-drift-aaa}
%\min_{t\in [t_0,t_0+T]} 
d(x(t),\cx^*) \ge \epsilon, ~t\in [t_0,t_0+T]
~~~~\mbox{implies}~~~~
\Phi(x(t_0+T))-\Phi(x(t_0)) \le -\delta.
\end{equation}
\end{lem}

{\em Proof.} If \eqn{eq-uniform-drift-aaa} would not hold, we would be able to construct 
a sequence of FSPs converging u.o.c. to an FSP such that 
$$
d(x(t),\cx^*)\ge \epsilon ~~\mbox{and}~~ \Phi(x(t)) = \Phi(x(0)), ~t\in [0,T].
$$
(Here we use the shift, continuity and uniform Lipschitz  properties of the 
family of FSPs, and the fact that $\cx$ is compact.)
This is not possible, because by Lemma~\ref{lem-drift-aaa} we must have
$(d/dt) \Phi(x(t)) < 0$ at every regular point in $[0,T]$.
$\Box$

As a corollary, we obtain the following
analog of Lemma~\ref{lem-fluid-convergence-closed}.

\begin{lem}
\label{lem-fluid-convergence-closed-aaa} 
Any FSP 
%(for the closed system under Greedy-D-AC) 
is such that
\beql{eq-fluid-convergence-closed-aaa}
d(x(t),\cx^*) \to 0.
\end{equation}
The convergence is uniform across all initial states $x(0)\in\cx$. 
\end{lem}

\subsection{Proof of of Theorem~\ref{th-fluid-stat-gt1-closed-aaa}.}

The rest of the proof of Theorem~\ref{th-fluid-stat-gt1-closed-aaa} is same as that of
Theorem~\ref{th-fluid-stat-gt1-closed}. 

\section{Some generalizations}
\label{sec-generalizations}

A number of generalizations of our results are not difficult to obtain.
We will discuss Theorems~\ref{th-fluid-stat-gt1-closed}
and \ref{th-fluid-stat-gt1-Greedy-DM} to be specific, but analogous generalizations
apply to Theorems~\ref{th-fluid-stat-gt1-closed-aaa} 
and \ref{th-fluid-stat-gt1-Greedy-DM-aaa}.

\subsection{A different procedure for placing arrivals.}

%Our main results -- let us talk about Theorems~\ref{th-fluid-stat-gt1-closed}
%and \ref{th-fluid-stat-gt1-Greedy-DM} to be specific --
Theorems~\ref{th-fluid-stat-gt1-closed} 
%(for Greedy-D in the closed system) 
and 
\ref{th-fluid-stat-gt1-Greedy-DM} 
%(for Greedy-DM in the original system),
 require that when a type $i$ customer (in Theorem~\ref{th-fluid-stat-gt1-closed})
or a type $i$ token (in Theorem~\ref{th-fluid-stat-gt1-Greedy-DM})
is assigned for service, it is placed along the edge $(k,i)$
minimizing the weight differential $\Delta_{ki}=\Delta_{ki}(X(t-))$.
%$W_{k+e_i}(X(t-))-I\{k\ne 0\}W_k(X(t-))$.
The procedure of choosing the edge to place a customer (or token) can be replaced
by the following one, which might be easier to implement in some scenarios.
We compare $\Delta_{k'i}$ for the edge $(k',i)$ along which a type $i$ departure
just occurred, to the $\Delta_{ki}$ for one edge, selected randomly as follows:
with probability $\epsilon\in (0,1)$ we select edge $(e_i,i)$; with probability
$1-\epsilon$ we pick a non-empty server uniformly at random and, if its 
configuration $\ell$ is such that $k=\ell+e_i \in \ck$, we select edge $(k,i)$. 
($\epsilon$ is a fixed parameter.) Now, if $\Delta_{ki}<\Delta_{k'i}$ for the
selected edge $(k,i)$ (if any), we place the customer (or token) along $(k,i)$;
otherwise, we place it ``back'' along $(k',i)$.
It is not difficult to see that the proofs of Theorems~\ref{th-fluid-stat-gt1-closed} 
and \ref{th-fluid-stat-gt1-Greedy-DM} still hold when Greedy-D and Greedy-DM algorithms,
respectively, are adjusted as described above.

The described alternative procedure generalizes the results in the sense that we can,
for example, use this procedure with a fixed probability $\delta \in [0,1]$
and use the the ``old'' procedure (picking the smallest differential $\Delta_{ki}$)
with probability $1-\delta$, and the results still hold.

\subsection{More general input processes and service time distributions.}
\label{sec-more-general}

Theorems~\ref{th-fluid-stat-gt1-closed} and 
\ref{th-fluid-stat-gt1-Greedy-DM} still hold for much more general 
input processes and
service time 
distributions (as opposed to Poisson and exponential, respectively).
For example,
a simple (but still far reaching) generalization is for the case
when, for each customer type $i$, 
the input process is renewal (i.i.d. interarrival times, with mean $1/(\lambda_i r)$
and finite variance) and
the service time distribution $G_i(\xi)$ (with mean $\int_0^{\infty} \xi dG(\xi)=1/\mu_i$) has
the ``hazard rate'' lower bounded by
$\mu_i^{min} \in (0,\mu_i]$: $dG(\xi)/[1-G(\xi)] \ge \mu_i^{min} d\xi, ~\forall \xi\ge 0$.
\iffalse
(That is, regardless of the elapsed service time,
the probability of a type $i$ customer completing service in a $dt$-long time interval
is at least $\mu_i^{min} dt + o(dt)$.)
\fi
In this case, we observe that the key conservation laws still hold for the fluid limit
of the stationary system: (a) the ``amount'' of type $i$ fluid is $\rho_i$ and remains
constant and (b) the total rate of (actual) type $i$ departures is $\lambda_i$ and remains
constant. In addition, 
say in the proof of
Theorem~\ref{th-fluid-stat-gt1-Greedy-DM}
to be specific, the corresponding FSPs are such that (actual) type $i$
departure 
rate from state $(k, \hat k)$ is lower bounded by $\mu_i^{min} \hat k_i x_{(k, \hat k)}(t)$.
(Of course, the FSPs need to be defined more
generally, to account for elapsed service times.) 
Given these properties, the entire argument goes through essentially as is. 
And, clearly, 
these properties hold under the input flow and service time assumptions
still far more general than in the simple case described above.
\iffalse
And we still have the key fact, which holds for general service time
distributions (with same means $1/\mu_i$):
the limit of stationary versions of the process
 is such that the total rate of $i$-departures at 
any time $t$ is exactly $\lambda_i$.
\fi

\section{Discussion}
\label{sec-discussion}

We have shown that (versions of) the Greedy algorithm are asymptotically
optimal 
%(for the closed and the original system, respectively) 
in the 
sense of minimizing the objective function $\sum_k X_k^{1+\alpha}$ with $\alpha>0$.
When $\alpha$ is small (but positive), the algorithms produce an approximation
of a solution minimizing the linear objective $\sum_k X_k$, i.e. the total number of occupied servers.
If $\sum_k X_k$ is the ``real'' underlying objective, the ``price'' we pay
by applying Greedy algorithm with small $\alpha>0$ is that the algorithm 
will keep non-zero amounts (``safety stocks'') of servers  in many
``unnecessary'' (from the point of view of linear objective)
 configurations $k$, including many -- potentially all -- non-maximal configurations in $\ck$.
What we gain for this ``price'' is the simplicity and agility 
of the algorithm. ``True'' minimization of the linear objective
$\sum_k X_k$ requires that a linear program is solved (via explicit offline or implicit dynamic approach),
so that the system is prevented from using ``unnecessary'' configurations $k$,
not employed in optimal LP solutions. 

The Greedy algorithm with $\alpha>0$ 
is asymptotically optimal as the average number $r$ of customers
in the system goes to infinity.
The fact that 
%in the limit
it  maintains safety stocks of many configurations, means in particular
that the algorithms' performance is close to optimal when the ratio
$r/|\ck|$ is sufficiently large, so that there is enough customers 
in the system to keep non-negligible safety stocks of servers in potentially all 
configurations. If the number $|\ck|$ of configurations is large,
then $r$ needs to be very large to achieve near-optimality.
The use of aggregate configurations in the special case of
vector-packing constraints 
%(which, as we have shown, preserves 
%the asymptotic optimality of Greedy algorithm), 
alleviates this scalability issue when the number $|\cq|$ 
of aggregate configurations is substantially smaller than $|\ck|$.

Finally, we note that the closed system, considered 
in Theorems~\ref{th-fluid-stat-gt1-closed} and \ref{th-fluid-stat-gt1-closed-aaa}, 
is not necessarily artificial.
For example, it models the scenario where VMs do not leave the system,
but can be moved (``migrated'') from one host to another.
% to improve certain objective.
In this case, a ``service completion'' is a time point when a VM migration can be attempted.
%If the migration process is governed by Greedy algorithm, then the minimization
%of $\sum_k X_k^{1+\alpha}$ is achieved (asymptotically, as the system size grows).

%\noindent
%{\bf Acknowledgments.} I would like to thank 

\iffalse

\newpage
\appendix
\begin{center}
{\Large\textbf{Appendix}}
\end{center}

\section{Proof of Lemma~\ref{lem-bound-tail}}
\label{prf-bound-tail}

\fi

\iffalse
\noindent
%Alexander L. Stolyar \\
Bell Labs, Alcatel-Lucent \\
600 Mountain Avenue, 2C-322 \\
Murray Hill, NJ 07974, USA \\
stolyar@research.bell-labs.com
\fi

\end{document}